\pdfoutput=1
 \documentclass[11pt]{article}

\usepackage{amsfonts}
\usepackage{graphicx}

\usepackage{cancel}

\usepackage{graphicx}
\usepackage[colorlinks=true, pdfstartview=FitV, linkcolor=blue, citecolor=blue,urlcolor=blue]{hyperref}

\usepackage{framed,color}
\definecolor{shadecolor}{rgb}{0.95, 0.95, 0.85}

\newcommand{\ba}{\begin{array}}
\newcommand{\ea}{\end{array}}

\newtheorem{corol}{Corollary}
\newcommand{\bcr}{\begin{corol}}
\newcommand{\ecr}{\end{corol}}

\newtheorem{hypo}{Assumption}
\newcommand{\bhypo}{\begin{hypo}}
\newcommand{\ehypo}{\end{hypo}}
\newtheorem{defi}{Definition}

\newtheorem{lemma}{Lemma}
\newcommand{\ble}{\begin{lemma}}
\newcommand{\ele}{\end{lemma}}

\newcommand{\bde}{\begin{defi}}
\newcommand{\ede}{\end{defi}}

\newtheorem{example}{Example}
\newcommand{\bex}{\begin{example}}
\newcommand{\eex}{\end{example}}

\newtheorem{prop}{Proposition}
\newcommand{\epr}{\end{prop}}
\newcommand{\bpr}{\begin{prop}}
\newtheorem{teo}{Theorem}
\newcommand{\bth}{\begin{teo}}
\newcommand{\eth}{\end{teo}}
\newtheorem{rema}{Remark}
\newcommand{\bre}{\begin{rema} \rm}
\newcommand{\ere}{ \end{rema}}

\newcommand{\bea}{\begin{eqnarray}}
\newcommand{\eea}{\end{eqnarray}}
\newcommand{\beas}{\begin{eqnarray*}}
\newcommand{\eeas}{\end{eqnarray*}}
\newcommand{\ee}{\end{equation}}
\newcommand{\be}{\begin{equation}}

\begin{document}

\centerline{\Large\bf  A Review on The Sixth Painlev\'e Equation
 } 

\vskip 0.3 cm
\centerline{\Large Davide Guzzetti}
\vskip 0.2 cm
\centerline{SISSA, Intenational School of Advanced Studies, 34136 Trieste, Italy}
\vskip 0.2 cm 
\centerline{guzzetti@sissa.it}



\begin{abstract}   
For the Painlev\'e 6 transcendents, we   provide a unitary description of the critical behaviours, the connection formulae,  their complete tabulation,  and the asymptotic distribution of poles close to a critical point.  
\vskip 0.2 cm
\end{abstract}

\vskip 0.5 cm 
\noindent
{\bf Note:} This paper is published in  Constr. Approx. 41 (2015), no. 3, 495–527. In the present version, a few typing mistakes have been corrected, together with other small details.

\vskip 0.5 cm 
\noindent
{\Large \bf Contents}

\vskip 0.3 cm
 The isomonodromy deformation method  provides a unitary description  of the critical behaviours of the solutions of the Painlev\'e 6 equation (PVI),  their connection formulae and the asymptotic distribution of poles close to a critical point.  The paper is a review of these results.  I have also included explanations  and connections among the above subjects which I have never had the chance to write in other papers.

{\bf -- Introduction:}  We introduce the Painlev\'e 6 equation, a non linear ODE which plays a central role in contemporary mathematics, and defines new non-linear special functions,  called {\it Painlev\'e Transcendents}.   We motivate the analysis of the {\it critical behaviours} of PVI transcendents close to the singular points of the equation ({\it critical points}), and the analysis of the {\it connection problem}. 

{\bf -- Section \ref{IDA}:}  We review the first general class of solution whose critical behaviours was discovered  by Jimbo, who solved the connection problem by means of the {\it method of monodromy preserving deformations}. Then, we  review the general scheme of this method. 

{\bf -- Section \ref{FCB}:}  We go deeper into the application of the  method of monodromy preserving deformations. First, we give further details on the procedure which provides the class of solutions due to Jimbo, then we review a matching procedure which produces critical behaviours not contained in the former class. We explain how the symmetries of PVI  allow to further enlarge the class of known transcendents. Finally, we collect all the critical behaviours known today in a {\it complete table}, in exactly the same way as the classical special functions were tabulated in the XIX and XX centuries.

{\bf --  Section \ref{Poles}:} We review the distribution of the {\it movable poles} of the Painlev\'e 6 transcendents close to a singular point of the equation. 

{\bf -- Section \ref{conver}:} We describe a method of local analysis, due to Shimomura, which provides the critical behaviours on  spiral shaped domains contained in the {\it universal covering of the punctured neighbourhood of a critical point}.   We also review the {\it Elliptic representation} of the transcendents, which yields the same results of Shimomura's approach. 

{\bf -- Section \ref{isoest}:}  We review how the isomonodromy deformation approach can be extended to the whole universal covering of the punctured neighborhood of a critical point, providing a comprehensive picture of the critical behaviours on this universal covering, including the position of the movable poles.

\section{Introduction}

\subsection{Background}
The Painlev\'e 6 equation, denoted PVI or $\hbox{PVI}_{\alpha \beta \gamma \delta}$, is the non linear differential equation  
$$
{d^2y \over dx^2}={1\over 2}\left[ 
{1\over y}+{1\over y-1}+{1\over y-x}
\right]
           \left({dy\over dx}\right)^2
-\left[
{1\over x}+{1\over x-1}+{1\over y-x}
\right]{dy \over dx}
$$
$$
+
{y(y-1)(y-x)\over x^2 (x-1)^2}
\left[
\alpha+\beta {x\over y^2} + \gamma {x-1\over (y-1)^2} +\delta
{x(x-1)\over (y-x)^2}
\right],~~~\alpha,\beta,\gamma,\delta\in{\bf C}
$$
This is the sixth and last of the non linear ODEs discovered by Painlev\'e \cite{pain} 
 and 
Gambier \cite{gamb}, 
who classified all the second order ordinary differential equations
of the type 
$$
         {d^2 y \over dx^2}= {\cal R}\left(x,y,{dy\over dx}\right)
$$
where ${\cal R}$ is rational in ${dy\over dx}$,  $x$ and
 $y$, such that the branch points and essential singularities depend only on the equation ({\it Painlev\'e Property}). 
The essential singularities and branch points  are called {\it critical points}; for PVI  
they are 0,1,$\infty$. The behaviour of a solution close to a critical point is called {\it critical 
behaviour}.
\footnote{This differs from the  terminology  of singularity theory, where a critical point is a zero of the first derivative of a function.}     
 A solution of the sixth Painlev\'e equation  
 can be analytically continued to a meromorphic function on the universal
   covering of ${\bf P}^1\backslash \{ 0,1 ,\infty \}$.

The six equations,  discovered at the beginning of the XX century, appeared to be irreducible to already known equations solvable in terms of classical functions. This fact has been rigorously proved only recently \cite{Umemura}, as we discuss below.   
In the last decades the Painlev\'e equations have emerged as one of the central objects in pure mathematics and mathematical physics, with   applications in a variety of problems, such as number theory, theory of analytic varieties (like Frobenius structures), random matrix theory, orthogonal polynomials, non linear evolutionary PDEs, combinatorial problems, etc. 
 The properties of the classical (linear) transcendental  functions have
been organised and tabulated in various classical handbooks.  A
comparable organisation and tabulation of the properties of the Painlev\'e functions is now needed. Today   we are able to write an essentially complete table of the critical behaviours, with full expansions, for the Painlev\'e 6 functions, and the corresponding connection formulae \cite{guz2012}.

\subsection{Solving PVI...}
 We need to solve the  non
 linear differential equation PVI. What does it mean that we know a solution? The functions which solve  PVI are generically "new
 transcendental functions",  called {\it Painlev\'e Transcendents}.     Umemura proved the  irreducibility of them to {\it
 classical functions}  \cite{Umemura} \cite{Umemura1} \cite{Umemura2}, namely functions  given in terms of a finite iteration of
 {\it permissible operations} applied to rational functions.  These operations are the derivation,  rational combination  (sum,
 product, quotient), algebraic combinations (the expression is a root of a polynomial whose coefficients are rational functions),
 contour integrals and quadratures, solution of a linear homogeneous differential equation whose coefficients are rational
 functions, 
        solution of an algebraic differential equation of the first order whose coefficients are rational functions, 
composition with abelian 
functions (the expression is  $\varphi(f_1(x),...,f_n(x))$, where $f_1$,...,$f_n$ are rational  functions, and $\varphi:{\bf C}^n/ \Gamma \to {\bf C}$ is meromorphic, $\Gamma$ is a lattice). The elementary transcendental functions are classical functions, because they are the algebraic functions, or a  function which is obtained from an algebraic function by integration (like the exponential, the trigonometric and hyperbolic functions), or the inverse of such an integral (like the logarithm, the elliptic functions, etc).
Umemura proved  that the general solution of a Painlev\'e equation is not a classical function.      
 H.Watanabe \cite{watanabe} applied the argument to PVI, and showed that a solution of PVI is either algebraic, or
solves a Riccati equation (one-parameter family of classical solutions), or it is not a classical  function.  All the algebraic solutions were classified in \cite{DM} when
  $\beta=\gamma=0$,
 $\delta={1\over 2}$, and then in \cite{Lisovyy} for the general PVI. 

Since a Painlev\'e transcendent is in general  not classical, we require the following minimal knowledge: 

{\bf ~i)} The knowledge of the  {\it explicit} critical behaviour (or the asymptotic expansion)  of a transcendent   at the critical points $x=0,1,\infty$.   
We symbolise the critical behaviour of {\it a branch} $y(x)$,  defined for  $-\pi<\hbox{arg} x <\pi$ and 
$-\pi<\hbox{arg} (1-x) <\pi$, as follows: 
 \be
\label{malamente1}
y(x)=y_u(x,{c_1}^u,{c_2}^u),~~~~~\hbox{ when } x\to u,~~~u\in\{0,1,\infty\},
\ee
 where  ${c_1}^u,{c_2}^u$ are the {\it integration constants}.

 {\bf ~ii)} The knowledge of the {\it explicit} connection formulae  
 among couples of integration
constants  at different critical points, as follows. 
A branch  $y(x)$ has a critical behaviour of type (\ref{malamente1}) at a critical point $x=u$, and another behaviour of type (\ref{malamente1}) at another critical point $x=v$, $u\neq v\in\{0,1,\infty\}$:
$$
y(x)=\left\{
\matrix{ y_u(x,c_1^{(u)},c_2^{(u)}),~~~x\to u
\cr 
\cr
 y_v(x,c_1^{(v)},c_2^{(v)}),~~~x\to v
}
\right.
$$
 The {\it connection problem} consists in the computation of the explicit formulae 
\be
\label{closd}
\left\{
\matrix{c_1^{(v)}=c_1^{(v)}(c_1^{(u)},c_2^{(u)})
\cr 
c_2^{(v)}=c_2^{(v)}(c_1^{(u)},c_2^{(u)})
}
\right. ,~~~\hbox{ and the inverse } ~ \left\{
\matrix{c_1^{(u)}=c_1^{(u)}(c_1^{(v)},c_2^{(v)})
\cr 
c_2^{(u)}=c_2^{(u)}(c_1^{(v)},c_2^{(v)})
}
\right. 
\ee
These are called the {\bf connection formulae in closed form}.

 Knowledge of {\bf i)} and {\bf ii)} is precisely what is meant by "solving"  a Painlev\'e equation in the review book \cite{Its}, page 8.  Above,  we used the word {\it
  explicit}. An explicit expression is 
a classical functions of its arguments. 
 Knowledge of {\bf i)} and {\bf ii)} allows to use Painlev\'e transcendents in applications. As already mentioned,  today   we
 are able to write an essentially complete table of the critical behaviours (with their full expansions) for the Painlev\'e 6
 transcendents, and the corresponding connection formulae \cite{guz2012}. We will review this in section \ref{tabulation}.

 The analysis of  the critical behaviours and asymptotic/formal expansions are  one side of the research on the Painlev\'e 6
 functions. Other perspectives, such as the classification of classical solutions \cite{watanabe}, rational solutions
 \cite{mazzoccoratio} and algebraic solutions \cite{DM} \cite{Lisovyy} will not be discussed here. 

\section{The Isomonodromy Deformation Approach}
\label{IDA}
The properties of some  Painlev\'e transcendents have been known since the beginning of the XX century, but the  first general result concerning the critical behaviour of a two parameter class of solutions of PVI
 is due to 
Jimbo in  \cite{Jimbo}.

\subsection{The Critical Behaviours in Jimbo's work (I)}
 A transcendent in the class obtained by  Jimbo    has critical behaviours \cite{Jimbo}: 
\be
\label{locintro}
y(x)=
\left\{
\matrix{
a_0 ~x^{1-\sigma_{0x}}\Bigl(1+O(x^\delta)\Bigr) ,& |x|<r
\cr 
\cr
1+a_1 ~(1-x)^{1-\sigma_{x1}}\Bigl(1+O((1-x)^{\delta})\Bigr) ,& |x-1|<r
\cr 
\cr
a_\infty ~x^{\sigma_{01}}\Bigl(1+O(x^{-\delta})\Bigr) ,&|1/x|<r
}
\right.
\ee
where $\delta$ is a small positive number,   $r>0$ is a sufficiently small radius,  $a_0,a_1,a_\infty \neq 0$ and 
 \be
0\leq\Re \sigma_{0x}<1,~~~0\leq\Re \sigma_{x1}<1,~~~0\leq\Re \sigma_{01}<1,~~~~~\sigma_{0x},\sigma_{x1},\sigma_{01}\neq 0
\label{RESTOacasa}
\ee 
are complex integration constats.   
The behaviours  (\ref{locintro}) hold when   $x$ converges
 to the critical points {\it inside a sector} with vertex
 on the corresponding critical point.  The angular width of the sector is arbitrary, and if increased the radius $r$ decreases. For  angular width $2\pi$, (\ref{locintro}) represents a branch of a Painlev\'e transcendent.  
The  two independent integration constants are any of  the three couples $(a_0,\sigma_{0x})$,  $(a_1,\sigma_{x1})$, $(a_\infty
,\sigma_{01})$.  The relation  among them is given by connection formulae (\ref{closd}), where the $c$'s represent the $a$'s and $\sigma$'s. 

The connection problem for (\ref{locintro}) was solved  in \cite{Jimbo},  when  
 $\alpha$, $\beta$, $\gamma$ $\delta$ are generic (we refer to  \cite{Jimbo} for a precise definition of 
 generic),     
  using the
{\it isomonodromy deformations theory}.  Jimbo considered an isomonodromic   $2\times 2$ 
Fuchsian system 
\be
   {d\Psi\over d\lambda}=A(x,\lambda)~\Psi,~~~~~
A(x,\lambda):=\left[ {A_0(x)\over \lambda}+{A_x(x) \over \lambda-x}+{A_1(x)
\over
\lambda-1}\right],~~~\lambda\in{\bf C}.
\label{SYSTEM}
\ee
with matrix coefficients  $A_i(x)$ ($i=0,x,1$).   
The isomonodromic deformation approach was  developed  in generality for the six Painlev\'e equations in \cite{JMU}  \cite{JM1}  \cite{JM2}, but for Fuchsian systems it goes back to R.Fuchs and  Schlesinger's work \cite{sch}.  The system (\ref{SYSTEM}) is by definition {\it isomonodromic} if there exists a fundamental $2\times2$ matrix solution $\Psi(\lambda,x)$ whose monodromy is independent of $x$. In other words, (\ref{SYSTEM}) must be one of the two systems of a   Lax Pair (note that this can be rephrased by  saying that we have an integrable structure).  The compatibility condition of the   Lax Pair is the  Schlesinger equations
$${dA_0\over dx}={[A_x,A_0]\over x},~~~{dA_1\over dx}={[A_1,A_x]\over 1-x},~~~{dA_x\over dx}={[A_0,A_x]\over x}+{[A_x,A_1]\over 1-x}
.   
$$
If the matrix coefficients of (\ref{SYSTEM}) satisfy special conditions,  the Schlesinger  equations are equivalent to the sixth Painlev\'e equation, as 
 established in \cite{JM1}.  These conditions are: 
\begin{equation}
 A_0+A_1+A_x = -{\theta_{\infty}\over 2}
 \sigma_3,~~\theta_\infty\neq 0.~~~~~
\hbox{ Eigenvalues}~( A_i) =\pm {1\over 2} \theta_i, ~~~i=0,1,x;
\label{caffe0}
\end{equation}
where the $\theta_\mu$'s, $\mu=0,x,1,\infty$, are defined in terms of the coefficients of  $\hbox{PVI}_{\alpha \beta \gamma \delta}$ as follows
 \begin{equation}
    \alpha= {1\over 2} (\theta_{\infty} -1)^2,
~~~-\beta={1\over 2} \theta_0^2, 
~~~ \gamma={1\over 2} \theta_1^2,
~~~ \left({1\over 2} -\delta \right)={1\over 2} \theta_x^2 ,~~~~~\theta_\infty\neq 0
\label{caffe1}
\end{equation}
Under the above conditions, the matrices $A_i(x)$ are given by explicit algebraic formulae in terms of $y(x)$. In particular, 
$$
A_{12}(x,\lambda)={g(x)(\lambda-y(x))\over \lambda(\lambda-1)(\lambda-x)},
$$
where $g(x)$ is a certain algebraic function of $x$. Therefore 
\be
\label{ipsilon}
y(x)={x(A_0)_{12}\over x[(A_0)_{12}+(A_1)_{12}]-(A_1)_{12}}
\ee
  The local behaviours (\ref{locintro}) are 
obtained  substituting into (\ref{ipsilon})   
 the critical  
 behaviours at $x=0$  of a class of  solutions   of the Schlesinger equations (see section \ref{JimboWork}).   The
 connection problem is then solved because the parameters $(a_0,\sigma_{0x})$,  $(a_1,\sigma_{x1})$ and  $(a_\infty,\sigma_{01})$  can be  expressed as functions of 
 the monodromy data of the Fuchsian system. We will come back to this point in sections \ref{mmm} and  \ref{JimboWork}.

\vskip 0.2 cm 
\noindent
{\bf Remark:} When $\Re \sigma_{ij}=0$,  there are three leading terms of the same order in $y(x)$. For example, at $x=0$:
$$
y(x)=a_0x^{1-\sigma_{0x}}+{4A^2\over a_0}x^{1+\sigma_{0x}}+Bx+O(x^{2-\sigma}),
$$
where $A=A(\sigma_{0x}^2)$ and $B=B(\sigma_{0x}^2)$ are
\be
\label{AB}
\left\{\matrix{
A^2:={(\sigma_{0x}^2-(\theta_0-\theta_x)^2)(\sigma_{0x}^2-(\theta_0+\theta_x)^2)\over 4\sigma_{0x}^2},\cr
\cr
B:={\theta_0^2-\theta_x^2+\sigma_{0x}^2\over 2\sigma_{0x}^2}.
}
\right.
\ee
We can rewrite the above as
\be
\label{hoshi0}
y(x)=x~\Bigl\{A\sin(i\sigma_{0x}\ln x+\phi)+B\Bigr\}+O(x^2),~~~a_0={A\over 2i}e^{i\phi}.
\ee
The last formula defines $\phi$ in terms of $a_0$, and the freedom $A\to -A$ corresponds to the freedom $\phi\to\phi+\pi$, while $a_0$ is fixed. 
 It may be convenient to write $\sigma_{0x}=-2i\nu$, $\nu \in {\bf R}$,  in such a way that
\be
\label{hoshi}
y(x)=x\left\{
{A}~\sin\bigl(2\nu\ln
x+\phi\bigr)
-{\theta_0^2-\theta_x^2-4\nu^2\over 8\nu^2}
\right\}+O(x^2),~~~~~x\to 0.
\ee
In the same way as above, if $\Re\sigma_{x1}=0$, then  (\ref{locintro}) at $x=1$ is replaced by 
\be
\label{HOS1}
y(x)=1-(1-x)\left\{A_1\sin\Bigl(i\sigma_{x1}\ln(1-x)+\phi_1\Bigr)+B_1\right\}+O((1-x)^2),
\ee
for suitable $A_1$ and $B_1$. 
If $\Re\sigma_{01}=0$, then  (\ref{locintro}) at $x=\infty$ is replaced (for suitable $A_\infty$ and $B_\infty$) by 
\be
\label{HOS2}
y(x)=A_{\infty}\sin\Bigl(i\sigma_{01}\ln\left({1\over x}\right)+\phi_\infty\Bigr) + B_\infty +O\Bigl({1\over x}\Bigr).
\ee

\subsection{Method of Monodromy Preserving Deformations}
\label{mmm}
 The paper \cite{Jimbo} opened the way to further research on PVI, based on the method of monodromy preserving deformations. The power of this method  is precisely that it allows to solve the connection problem. Its general scheme  is as follows. 

\begin{figure}
\centerline{\includegraphics[width=0.5\textwidth]{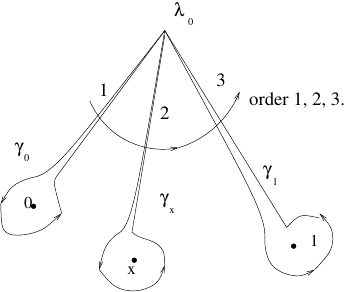}}
\caption{The ordered basis of loops}
\label{figure1}
\end{figure}

First, we define the monodromy data of the system (\ref{SYSTEM}). 
In  the ``$\lambda$-plane'' 
${\bf C}\backslash\{0,x,1\}$ we fix a base point $\lambda_0$ 
and  three loops, which are numbered in order 1, 2, 3 according to a 
counter-clockwise order referred to $\lambda_0$.  We choose $0,x,1$ to be the order $1,2,3$. 
We denote the loops by $\gamma_0$, $\gamma_x$, $\gamma_1$. See figure \ref{figure1}.   
  The loop at infinity will be $\gamma_\infty=\gamma_0\gamma_x\gamma_1$. 
  When $\lambda$ goes around a small loop around $\lambda=i$, $i=0,x,1$,  the  fundamental solution transforms like $\Psi \mapsto \Psi M_i$, where $M_0$, $M_x$, $M_1$ are the monodromy matrices w.r.t. the base of loops.  
Let $(\theta_0, \theta_1, \theta_x,\theta_\infty)\in{\bf C}^4$ be fixed by PVI, up to the equivalence $\theta_k\mapsto -\theta_k$, $k=0,x,1$, and $\theta_\infty\mapsto 2-\theta_\infty$. Denote $\sim$ the equivalence and let  the quotient be
$$
\Theta:=\{(\theta_0, \theta_1, \theta_x, \theta_\infty)\in{\bf C}^4~|~\theta_\infty\neq 0\}/ \sim.
$$  
  Let $M_{\infty}:=M_1M_xM_0$ be the monodromy at $\lambda=\infty$,  and consider  the set of triples of (invertible) monodromy matrices with unit determinant, defined up to conjugation $M_i\mapsto CM_i C^{-1}$ ($i=0,x,1$) by an invertible matrix $C$, namely  
$$ 
  M:=\{(M_0,M_x,M_1)\in\hbox{SL}(2,{\bf C})~| ~\hbox{Tr}M_\mu= 
2\cos\pi\theta_\mu,~\mu=0,1,x,\infty\}/\hbox{conjugation}
$$

\vskip 0.2 cm 
\noindent
{\bf Definition}: The {\it  monodromy data} of the class of Fuchsian systems (\ref{SYSTEM}), with the basis of loops ordered as  figure \ref{figure1}, is the set ${\cal M}:=\Theta\cup M$. 
\vskip 0.2 cm 
\noindent

When we fix  branch cuts in the $x$ plane, for example  $-\pi<\arg x<\pi$ and $-\pi<\arg (1-x) <\pi$, then   
to every branch $y(x)$, a  system (\ref{SYSTEM}) is associated, and so is a point in ${\cal M}$. Conversely, to a point in ${\cal M}$ a system or a family of systems (\ref{SYSTEM}) is associated through a Riemann-Hilbert problem \cite{AB}, and so is  either one branch $y(x)$ or a family of branches $y(x)$. Let  
\be
\label{monmap}
 f:\{y(x)\hbox{ branch}\} \rightarrow {\cal M}
\ee
be the map from the set $\{y(x)\hbox{ branch}\}$  of all the branches  of all equations $\hbox{PVI}_{\alpha,\beta,\gamma,\delta}$, $(\alpha,\beta,\gamma,\delta)\in{\bf C}^4$,  onto ${\cal M}$, associating to a branch the corresponding monodromy data.  
This map gives a way of classifying the branches of PVI transcendents and their critical behaviours in terms of monodromy data. It is through this map that we can in principle answer  the question wether the class of critical behaviours known today is complete.   For if the image through $f$ of this class is the whole of ${\cal M}$, then we are sure that the class contains all the possible critical behaviours. $f$ is generically injective, according to the following

\bpr 
\label{injec} Let the order of loops be fixed. The map (\ref{monmap})  is injective (one-to-one) when restricted to $f^{-1}(\Theta\cup \{(M_0,M_x,M_1)\in M ~|~M_\mu\neq I,~~\forall\mu=0,x,1,\infty\})$ \cite{D1}.
\epr

 \vskip 0.2 cm 
Let (\ref{malamente1}) be the the critical behaviour for $x\to u$, $u\in\{0,1,\infty\}$, 
 depending on two integration constants ${c_1}^u,{c_2}^u$. The map $f$ can be  made explicit because the monodormy matrices    
can be  computed as explicit functions of  ${c_1}^u,{c_2}^u$.   Conversely,   ${c_1}^u,{c_2}^u$ can be  expressed in terms of the 
associated monodromy data, provided that  $f$ is restricted as in Proposition \ref{injec}, so that it is injective.  
In order to make $f$ explicit, we need parameters (coordinates) to identify a point of ${\cal M}$.  Let 
$$
p_{ij}:=\hbox{Tr}(M_iM_j), ~~~j=0,x,1;~~~~~~ p_\mu:=\hbox{Tr}M_\mu=2\cos\pi\theta_\mu, ~~~\mu=0,x,1,\infty. 
$$
Observe that $p_{ij}=p_{ji}$. These are seven invariant functions (w.r.t. conjugation and $\sim$) defined on ${\cal M}$.  They satisfy  the  relation 
$$
p_{0x}^2+p_{01}^2+p_{x1}^2+p_{0x}p_{01}p_{x1}-(p_0p_x+p_1p_\infty)p_{0x}-(p_0p_1+
p_xp_\infty)p_{01}-(p_xp_1+p_0p_\infty)p_{x1}+
$$
\be
\label{FR}
+p_0^2+p_1^2+p_x^2+p_\infty^2+p_0p_xp_1p_\infty -4 =0
\ee
 in agreement with the dimension of ${\cal M}$.  This relation appeared in the 1897 book of Fricke and Klein   \cite{FK} and was rediscovered by Jimbo in \cite{Jimbo} (it follows from the trace of the relation $M_1M_xM_0=M_\infty$). The geometry of ${\cal M}$, where (\ref{FR}) is an affine cubic,  was  studied in \cite{Iwa}. 
   The $p_{ij}$ are the local coordinates on ${\cal M}$ we are looking for, according to the following proposition:  
\bpr
\label{Jimcoo}
    $p_0,p_x,p_1,p_\infty$; $p_{0x}, p_{x1}, p_{1x}$     are   coordinates on the subset of ${\cal M}$ where the group generated by  $M_0$, $M_x$, $M_1$ is irreducible  \cite{Iwa}.
\epr

We come back to the problem of the explicit form of $f$. Suppose that  Propositions \ref{injec} and \ref{Jimcoo} hold for a given $\hbox{PVI}_{\alpha,\beta,\gamma,\delta}$. The $p_\mu$'s  $\mu=0,x,1,\infty$ are fixed by $\alpha,\beta,\gamma,\delta$ and only two parameters $p_{ij}$'s are independent. They may be thought as the integration constants for the associated branch $y(x)$  with critical behaviour (\ref{malamente1}).  The asymptotic techniques of the method of monodromy preserving deformations provide explicit {\it parametric formulae}  
\be
\label{POA}
\left\{
\matrix{ {c_1}^u={c_1}^u(p_0,p_x,p_1,p_\infty,p_{0x},p_{x1},p_{01})
\cr 
{c_2}^u={c_2}^u(p_0,p_x,p_1,p_\infty,p_{0x},p_{x1},p_{01})
}
\right.
\ee
and the inverse formulae 
\be
\label{POA1}
\left\{
\matrix{
p_{0x}=p_{0x}({c_1}^u,{c_2}^u,p_0,p_x,p_1,p_\infty) 
\cr 
p_{x1}=p_{x1}({c_1}^u,{c_2}^u,p_0,p_x,p_1,p_\infty) 
\cr 
p_{01}=p_{01}({c_1}^u,{c_2}^u,p_0,p_x,p_1,p_\infty) 
}
\right.
\ee
Explicit means that the formulae are classical functions of their arguments, as we have already explained.  
 Three pairs of different  parametric formulae of type (\ref{POA}) can be written at $x=0,1$ and $\infty$ respectively,   
in terms of the same monodromy data, namely for the same branch $y(x)$. Conversely, the monodromy data associated to  a given 
$y(x)$ can be written as in (\ref{POA1})  in three ways, namely in terms of the three couples of integration constants  at $x=0$, 
$1$ and $\infty$ respectively. Thus, we say that the formulae  (\ref{POA}) and   (\ref{POA1})  are  {\bf parametric connection formulae}. The {\bf connection formulae in closed form} (\ref{closd}) are obtained form the parametric form, combining (\ref{POA}) and (\ref{POA1}) at two critical points $x=u$ and $v$.

 An example of the above parametric formulae is the  relation, established in \cite{Jimbo},  among the {\it monodromy exponents} $\sigma_{0x}$, $\sigma_{x1}$ and $\sigma_{01}$  of (\ref{locintro}) and the $p_{ij}$'s:
\be\label{STA}
2\cos\pi\sigma_{0x}=p_{0x},~~~~~2\cos\pi\sigma_{x1}=p_{x1},~~~~~2\cos\pi\sigma_{01}=p_{01}
\ee
It is clear from the above that the restriction (\ref{RESTOacasa}) implies that the image through $f$ of the branches with behaviour (\ref{locintro}) and (\ref{hoshi0}) (and similar behaviours like (\ref{hoshi0})  at $x=1,\infty$)  is contained in a wide  subspace of ${\cal M}$, characterised by  
$$ p_{0x},p_{x1},p_{01}\not\in (-\infty,-2]\cup\{2\}.$$ 
In particular, for solutions  (\ref{hoshi0}) we have  $p_{0x}>2$, while $p_{x1}>2$, $p_{01}>2$ for  (\ref{HOS1}) and (\ref{HOS2})  respectively. 
 Similalrly, there are explicit formulae 
\be
\label{STA1}
a_i=a_i(p_0,p_x,p_1,p_\infty,p_{0x},p_{x1},p_{01})\equiv 
a_i(\theta_0,\theta_x,\theta_1,\theta_\infty,\sigma,p_{x1},p_{01}),~~~i=0,1,\infty,
\ee
which the reader can find in \cite{Jimbo}, or more  explicitly  in \cite{Boalch} (they are collected in \cite{guz2012} as well). 

\section{Finding the Critical Behaviours}
\label{FCB}

In section \ref{IDA}, we have described a wide class of critical behaviours (\ref{locintro}) and (\ref{hoshi0}). In this section, we explain in more details how they were derived. Then, we discuss the matching method, which allows to find further critical behaviours. We also review the symmetries of PVI, which, applied to the critical behaviours obtained by Jimbo method and the matching method, produce further critical behaviours. As a final result, we obtain essentially all the possible critical behaviours, with their full expansion. They can be tabulated  like the classical special functions, as we will see in section \ref{tabulation}.  

\subsection{Preamble: Symmetries}
 It is possible to study critical behaviours and parametric connection formulae at one critical point, say $x=0$, and transfer the results to the others. This is because PVI admits symmetry transformations, called 
{Okamoto's bi-rational transformations}, introduced in \cite{Okamoto}, of the form 
\be
\label{OKAvino}
y^{\prime}(x)={P\left(x,y(x),{dy(x)\over dx}\right)\over Q\left(x,y(x),{dy(x)\over dx}\right)}, ~~~x^{\prime}={p(x)\over q(x)},~~~~~
(\theta_0,\theta_x,\theta_1,\theta_\infty)\mapsto 
(\theta_0^{\prime},\theta_x^{\prime},\theta_1^{\prime},\theta_\infty^{\prime})
\ee
such that $y(x)$ satisfies (PVI) with coefficients $\theta_0,\theta_x,\theta_1,\theta_\infty$ and variable $x$,  if and only if $y^{\prime}(x^{\prime})$ satisfies (PVI) with coefficients $\theta_0^{\prime},\theta_x^{\prime},\theta_1^{\prime},\theta_\infty^{\prime}$ and variable $x^\prime$. The functions 
 $P,Q$ are polynomials; $p,q$ are linear; the transformation of the $\theta_\mu$'s
  is an element of a linear representation of the Weyl Group  of the root system $D_4$, or the group of  Shift $l_j: v=(v_1,v_2,v_3,v_4)\mapsto 
v+e_j$, $j=1,2,3,4$, where  $e_1=(1,0,0,0)$, ..., $e_4=(0,0,0,1)$, or  the Permutation group. The latter acts as follows:  
\be
\theta_0^{\prime}=\theta_1,~~\theta_x^{\prime}=\theta_x,~~\theta_1^{\prime}=\theta_0,~~\theta_\infty^{\prime}=\theta_\infty;~~~y^{\prime}(x^\prime)=1-y(x),~~~x^\prime=1-x. 
\label{onara}
\ee
 \be
\label{onara1}
\theta_x^\prime=\theta_1,~~\theta_1^\prime=\theta_x;~~~~~~~
\theta_0^\prime=\theta_0,~~\theta_\infty^\prime=\theta_\infty;~~~~~~~~~~
y^\prime(x^\prime)={1\over x}y(x),~~~x^\prime={1\over x}.
\ee
\be
\label{sym2}
\theta_0^\prime=\theta_\infty-1,~~\theta_x^\prime=\theta_1,~~\theta_1^\prime=\theta_x,~~\theta_\infty^\prime=\theta_0+1;~~~~~~~y^\prime(x^\prime)={x\over y(x)},~~~x=x^\prime.
\ee 
 (\ref{onara}) and (\ref{onara1})   produce critical behaviours at $x=1,\infty$ respectively, from behaviours at $x=0$. 
The transformation (\ref{sym2})  generates new behaviours from known ones at a given critical point. 

 The effect of the bi-rational transformation on the  monodromy data $p_{ij}$ is described in \cite{Jimbo}, \cite{DM}, \cite{D1}, and more  generally in \cite{DM11} (see also \cite{Marta}). In particular, we have: 
\be
\label{ricordo11}
\hbox{For (\ref{onara}): }~~ p_{01}^\prime
=-p_{01}-p_{0x}p_{x1}+p_\infty p_x+p_1p_0 ,~~~
p_{0x}^\prime=p_{x1}~~~
p_{x1}^\prime=p_{0x}.
\ee
\be
\label{ricordo22}
\hbox{For (\ref{onara1}) :}~~ 
p_{0x}^\prime =
 -p_{01}-p_{0x}p_{x1} +p_\infty p_x +p_0p_1 ,~~~
p_{01}^\prime=p_{0x},~~~
p_{1x}^\prime=p_{1x}.
\ee
\be
\label{ricordo3}
\hbox{For  (\ref{sym2}):}~~~~~~~~~~~~~~~~~~~~~~~~~~~~~~~
p_{0x}^\prime=-p_{0x},~~~p_{01}^\prime=-p_{01},~~~p_{x1}^\prime=p_{x1}.
\ee

\subsection{The Critical  Behaviours  in Jimbo's work (II)}
\label{JimboWork}

The behaviours (\ref{locintro}) (or  (\ref{hoshi0}),  (\ref{HOS1}) and (\ref{HOS2}))  follow  from   a result on  a class of solutions  of the Schlesinger equations, established in Lemma 2.4.8 (page 262) of \cite{SMJ}, stating the following. 

Let  $A^0_0$, $A^0_x$, $A_1^0$ be constant matrices satisfying (\ref{caffe0}). Let  $\Lambda:=A^0_0+A^0_x$ and let its eigenvalues $\pm\sigma_{0x}/2$ be given (they will be precisely the monodromy exponents in (\ref{locintro})). Suppose that $|\Re \sigma_{0x}|<1$ and $\Lambda\neq 0$.  Then, for any  $\vartheta>0$ there exists a small $r=r(\vartheta)$ such that the Schlesinger equations admit a {\it unique} solution defined for $0<|x|<r$, $|\arg(x)|<\vartheta$, with the properties 
\be
\label{Lambda}
A_1(x)=A^0_1+O(|x|^{\delta}),~~~x^{-\Lambda}A_0(x)~x^\Lambda=A^0_0+O(|x|^{\delta}),
~~~x^{-\Lambda}A_x(x)~x^\Lambda=A^0_x+O(|x|^{\delta})
\ee
where $\delta\leq 1-|\Re\sigma_{0x}|$ is a positive small number. Note that the following limit exists:
$$
\Lambda=\lim_{x\to 0} (A_0(x)+A_x(x))\neq 0,
$$
It is not difficult to find the general parametrization of  $A^0_0$, $A^0_x$, $A_1^0$ satisfying (\ref{caffe0}) and $A^0_0+A^0_x=\Lambda$. An extra parameter $a_0\neq 0$ will appear, together with $\sigma_{0x}$ . 
 Then, the behaviour (\ref{locintro}) at $x=0$ (or (\ref{hoshi0}), if $\Re \sigma_{0x}=0$),  is obtained by the substitution of
 the critical behaviours (\ref{Lambda})  into  (\ref{ipsilon}).  The behaviours at $x=1,\infty$ can then be computed by means of 
  (\ref{onara}) and (\ref{onara1}).  Note that one can substitute the condition $|\Re\sigma_{0x}|<1$ with the equivalent
 $$0\leq \Re\sigma_{0x}<1.$$ 

Jimbo obtained the connection formulae  (\ref{POA}) and (\ref{POA1}) (the role of  $(c_1^{0},c_2^0)$ is played by $(\sigma_{0x},a_0)$), by  reducing the system  ({\ref{SYSTEM}) to two simpler sistems, whose monodromy matrices are exactly computable.  He proved that  the following limits, constructed from an isomonodromic  fundamental solution $\Psi(\lambda,x)$  of (\ref{SYSTEM}), exist 
$$
\Psi_{OUT}(\lambda):=\lim_{x\to 0} \Psi (\lambda, x),~~~~~\Psi_{0}(\mu):=\lim_{x\to 0} x^{-\Lambda}\Psi(x\mu,x),~~~~~~\lambda=x\mu
$$
 and satisfy
\be
\label{fuchsianSYSTEM}
{d\Psi_{OUT}\over d \lambda}=
\left[
{\Lambda\over \lambda}+ {A^0_1\over \lambda-1}
\right]~\Psi_{OUT},~~~
{d \Psi_{0} \over d\mu}
=
\left[
{A^0_0\over \mu} + {A^0_x \over \mu -1}
\right]~
\Psi_{0}.
\ee 
 The notation $\Psi_{0}$ and $\Psi_{OUT}$ is not  in the original paper \cite{Jimbo};  we will explain the reason for it below.  
It follows from the isomonodromy property that  $M_1$   coincides with the monodromy $M_1^{OUT}$ of $ \Psi_{OUT}$, while $M_0$ and $M_x$ coincide with the monodromies at $\mu=0$ and $\mu=1$  of $\Psi_{IN}$ respectively. We call these monodromies $M_0^{IN}$ and $M_x^{IN}$.  Systems (\ref{fuchsianSYSTEM})   are equivalent to Gauss
hyper-geometric equations and the monodromy matrices $M_1^{OUT}$, $M_0^{IN}$ and $M_x^{IN}$ can be computed from the connection formulae for the hyper-geometric functions. Their matrix elements  contain $\sigma_{0x}$ and  $a_0$ as parameters. Then,  (\ref{POA}) and (\ref{POA1}) for (\ref{locintro}) follow by taking the traces of the products $M_iM_j$, $i,j\in{0,x,1}$. We refer to  the original works \cite{Jimbo} and \cite{Boalch} (and the collection/tabulation  \cite{guz2012}) for the resulting formulae  (\ref{STA}) and  (\ref{STA1}). 

Observe that for $x\to 0$ one can also expect that the system (\ref{SYSTEM}) can be approximated by a system 
$$
{d \Psi_{IN} \over d\lambda}
\cong 
\left[
{A_0(x)\over \lambda} + {A_x(x) \over \lambda -x}
\right]~
\Psi_{IN},
$$
when $x\to 0$ and also $\lambda\to 0$, so that  $A_1/(\lambda-1)$ is neglected. Note that in case of  behavior (\ref{Lambda}), the matrixes $A_0(x)$ and $A_x(x)$ diverge for $x\to 0$.
\subsection{Matching}
\label{Matching}

 The critical behaviours (\ref{locintro}) and (\ref{hoshi0}) are  based on the restriction (\ref{RESTOacasa}), which is a restriction on $\Lambda$. The question arises  whether the limit procedure of section \ref{JimboWork} can be generalised  when $\Lambda=\pmatrix{0&1\cr0&0}$ and when 
$$
\lim_{x\to 0} (A_0(x)+A_x(x))=0
$$ 
We also would like to study cases where $A_0(x)$ and  $A_x(x)$, taken separately, diverge faster than a power $x^{-|\sigma|/2}$, $0<|\sigma|<1$, for example, cases when $xA_0(x)\to $ constant matrix, for $x\to 0$.  These  behaviours substituted into (\ref{ipsilon}) will produce critical behaviours different from (\ref{locintro}). 

 The above question, together with the search for a constructive procedure which  produces critical behaviours of  $A_0(x)$, $A_x(x)$ and $A_1(x)$  with minimal assumptions, was the starting point for the {\it matching procedure} developed in \cite{D2}. Originally, such procedure   
was suggested by Its and Novokshenov in \cite{IN}, 
 for the second and third Painlev\'e equations. 
 This approach was further developed and
used by Kapaev, Kitaev,  Andreev, and Vartanian (see for example \cite{KitaevAndreev}). We finally implemented it  for  PVI in \cite{D2} (I thank A. Kitaev for introducing me to the method).

 We briefly explain the construction. Let us divide the $\lambda$-plane into
two domains. The ``outside'' domain  is defined for $
|\lambda|\geq |x|^{\delta_{OUT}}$, $\delta_{OUT}>0$. In this domain,  (\ref{SYSTEM}) can be approximated by: 
\be
{d\Psi_{OUT}\over d \lambda}=
\left[
{A_0+A_x\over \lambda}+{xA_x\over \lambda^2}~\sum_{n=0}^{N_{OUT}}\left({x\over \lambda}\right)^n+ {A_1\over \lambda-1}
\right]~\Psi_{OUT},
\label{nonfuchsianSYSTEMOUT}
\ee
The ``inside'' domain  is defined for $\lambda$ comparable with $x$, namely $
|\lambda|\leq |x|^{\delta_{IN}}$, $\delta_{IN}>0$, and ({\ref{SYSTEM}) can be approximated by: 
 \be
{d \Psi_{IN} \over d\lambda}
=
\left[
{A_0\over \lambda} + {A_x \over \lambda -x} - A_1 \sum_{n=0}^{N_{IN}} 
\lambda^n
\right]~
\Psi_{IN},
\label{nonfuchsianSYSTEMIN}
\ee
where $N_{IN}$, $N_{OUT}$ are suitable integers.
 The leading term of $y(x)$  is  obtained  by requiring that the fundamental matrix 
 solutions $\Psi_{OUT}(\lambda,x)$, $\Psi_{IN}(\lambda,x)$ 
match in the region of overlap, provided this is not empty: 
\be\Psi_{OUT}(\lambda,x)
\sim 
\Psi_{IN}(\lambda,x), ~~~~~
|x|^{\delta_{OUT}} \leq |\lambda|\leq |x|^{\delta_{IN}},~~~x\to 0  
\label{overlapINOUT}
\ee  
This relation is to be intended in the sense that the leading terms 
of the local behaviour of $\Psi_{OUT}$ and $\Psi_{IN}$ for $x\to
0$ must be equal (note: $\delta_{IN}\leq\delta_{OUT}$).  In this procedure,  also the reduced systems (\ref{nonfuchsianSYSTEMOUT}) and (\ref{nonfuchsianSYSTEMIN}) are isomonodromic.   
The  matching condition (\ref{overlapINOUT}) produces  the $x$-leading terms for $x\to 0$  of $\Psi_{OUT}(\lambda,x)$ and $\Psi_{IN}(\lambda,x)$, and then of $A_0(x)$, $A_1(x)$ and $A_x(x)$ (recall that $A=d\Psi/d\lambda \cdot \Psi^{-1}$). Finally,  it produces the leading term(s) of $y(x)$ in (\ref{ipsilon}).  

The  Fuchsian systems (\ref{fuchsianSYSTEM}) represent the simplest case, occurring   when $0\leq \Re \sigma_{0x}<1$, but the procedure in general involves non Fuchsian reductions. As a result,  new critical behaviours, different from (\ref{locintro}) and (\ref{hoshi0}),  are produced through  (\ref{ipsilon}). For example: 

-- The (basic) logarithmic behaviours when $\sigma_{0x} = 0$: 
\be
\label{intrlog1}
y(x)=
x\left[
{\theta_x^2-\theta_0^2\over 4} 
\bigl(
\ln x +a
\bigr)^2 +{\theta_0^2\over \theta_0^2-\theta_x^2}
\right]+O(x^2\ln^3x)
,~~~~~\hbox{ if } \theta_0^2\neq\theta_x^2,
\ee
\be
\label{intrlog2}
y(x)= x(a~\pm~\theta_0~\ln x)+O(x^2\ln^2x),~~~~~~~\hbox{ if } \theta_0^2=\theta_x^2,
\ee
where $a$ is an integration constant. This solution is discussed in \cite{Jimbo} as well, as a limit of (\ref{locintro}) for $\sigma_{0x}\to 0$. 

-- The (basic) Taylor expansions: we will write them in the table of section  \ref{tabulation}. 
%
%
%

\vskip 0.2 cm
 In order to compute the parametric connection formulae , we compute the  monodromy data of ({\ref{SYSTEM})  as follows.  {\it 
Once the matching $\Psi_{OUT}\leftrightarrow \Psi_{IN}$ in (\ref{overlapINOUT}) has been completed}, we  match $\Psi_{OUT}$ with a fundamental solution $\Psi$ of ({\ref{SYSTEM}) at $\lambda =\infty$ and $\lambda=1$. Namely $\Psi_{OUT}$  and $\Psi$ match in the domain $|\lambda|\geq |x|^{\delta_{OUT}}$, which does not contain $\lambda=0$ and $\lambda=x$. 
We also  match $\Psi_{IN}$ with {\it the same} $\Psi$ in another region of the $\lambda$-plane, around  $\lambda=0$ and $x$, namely  in the domain $|\lambda|\leq |x|^{\delta_{IN}}$, which does not contain $\lambda=1$ and $\lambda=\infty$.

   If this matching is realised
\footnote{Note that $\Psi_{IN}$ does not match with $\Psi$ at $\lambda=1$ and $\infty$, and $\Psi_{OUT}$ does not match with $\Psi$ at $\lambda=0$ and $x$. 
},   $M_1$  coincides with $M_1^{OUT}$ of $ \Psi_{OUT}$, while $M_0$ and $M_x$ coincide with $M_0^{IN}$ and $M_x^{IN}$ of $\Psi_{IN}$. Their entries contain the integration constants of the critical behaviour of $y(x)$. Thus,  (\ref{POA}) and (\ref{POA1})  follow by taking the traces of the products $M_iM_j$, $i,j\in{0,x,1}$.  Since we need  to  compute the monodromy matrices $M_1^{OUT}$, $M_0^{IN}$ and $M_x^{IN}$ exactly, we should be able to solve the reduced  "IN" and "OUT" systems in terms of linear special functions, or find that they are the linear systems associated to other Painlev\'e equations for which we already know the monodromy.  In \cite{D2} and \cite{D1} the monodromy matrices associated to  the basic Taylor and logarithmic solutions are computed. The Taylor solutions and  the  associated  monodromy are computed also in \cite{Kaneko}.

\subsection{Other behaviours generated by symmetries} 

 The birational transformations (\ref{onara}) and (\ref{onara1}) generate the critical behaviours at $x=1$ and $\infty$. For example, when applied to 
(\ref{intrlog1}), they provide the behaviours:
$$
y(x)\sim
 1 -(1-x) \left\{{\theta_1^2\over \theta_1^2-\theta_x^2} +
{\theta_x^2-\theta_1^2\over 4}
\bigl(
\ln (1-x) +a_1
\bigr)^2
\right\},~~~x\to 1,
$$
$$
y(x)\sim
 {\theta_0^2\over \theta_0^2-\theta_1^2}+
{\theta_1^2-\theta_0^2\over 4} 
\left[
\ln {1\over x} +a_\infty
\right]^2,~~~x\to\infty.
$$
 
 The symmetry (\ref{sym2}) produces new behaviours at a given critical point starting from known ones. For example, when applied to (\ref{intrlog1})  at $x=0$, it provides the branch with the behaviour
$$
y(x)={4\over (\theta_1^2-(\theta_\infty-1)^2)~\ln^2 x}\left[1-{2a\over \ln x} +\left({1\over \ln^2x}\right)\right],~~~x\to 0
$$
When applied to (\ref{hoshi}), it generates a branch with behaviour
\be
\label{invos1}
y(x)={1\over {\cal A}\sin(2\nu\ln x+\phi)+{\cal B}+O(x)}
, ~~~~~x\to 0
\ee
where  $\nu\in {\bf R}\backslash\{0\}$ and $\phi\in{\bf C}$ are integration constats, 
$${\cal A}= \left[
{(\theta_\infty-1)^2\over 4\nu^2}+\left(
{(\theta_\infty-1)^2-\theta_1^2-4\nu^2\over 8\nu^2}
\right)^2
\right]^{1\over 2},~~~~~{\cal B}= -{(\theta_\infty-1)^2-\theta_1-4\nu^2\over 8\nu^2},
$$ and $p_{0x}=2\cos(\pi(1+2i\nu))<-2$. 
There are similar behaviours at $x=1$ and $\infty$, namely 
\be
\label{invos2}
 y(x)=1-{1\over {\cal A}_1\sin(2\nu_1\ln(1-x)+\phi_1)+{\cal B}_1+O(1-x)},~~~x\to 1; ~~~p_{x1}<-2,
\ee
\be
\label{invos3}
y(x)={x\over {\cal A}_\infty \sin(2\nu_\infty \ln(1/x)+\phi_\infty)+{\cal B}_\infty +O(1/x)},~~~x\to\infty;~~~p_{01}<-2.
\ee
with suitable ${\cal A}$ and ${\cal B}$.

\subsection{Tabulation}
\label{tabulation} 
 The critical behaviours admit  full expansions (convergent, asymptotic or formal). 
 For example,  the  procedure to compute the   expansion of  (\ref{locintro}) at $x=0$ is given  in \cite{guz2010}. It involves a recursive computation, based on the substitution into PVI,  which increases in complication with the order of the expansion.  The Taylor expansions 
   can be obtained by a recursive computation based on the substitution into PVI. The same can be done starting from the behaviours (\ref{intrlog1}) and (\ref{intrlog2}), to generate  two full "series expansions"  at $x=0$ whose coefficients are polynomials of $\log x$.  

It is to be remarked that the closed form of the coefficients  {\it at all orders}  of  the expansion of the $\tau$ function associated to  (\ref{locintro}) when $x\to 0$ was first given  in \cite{Lisovyy1}, in terms of conformal blocks.  This resut is reviewed in this Special Issue of CA, in paper \cite{Lisovyy2}. 

 The  birational transformations  associated to the Weyl group,  the Shifts and Permutations  applied to  the full expansion of 
the solutions obtained in subsections \ref{JimboWork}  and \ref{Matching}, generate a variety of critical behaviours, whose image 
through $f$  is essentially the whole of  ${\cal M}$ \cite{guz2012}.
 The systematic tabulation of these behaviours is possible. In the paper 
\cite{guz2012}, we  provided the tables of the critical behaviours with full expansion at $x=0,1$ and $\infty$, and  the parametric connection 
formulae. 
Below, we reproduce the table of \cite{guz2012}, only at  $x=0$, and in a simplified form. 
The branches in the table may be classified according to their behaviour: 

\vskip 0.2 cm 
- Complex  power behaviours: they are expanded in powers of $x^{n+m\lambda}$, for some $n,m\in{\bf Z}$ and $\lambda\in{\bf C}$. In this case,   $|y(x)|$ may vanish, converge to a constant or diverge when $x\to 0$.  They include (\ref{locintro}) and 
(\ref{hoshi0}).

\vskip 0.2 cm 
- Inverse oscillatory behaviours: $y(x)$ oscillates without vanishing when $x\to 0$, and may have poles in a sector centred at $x=0$ (these poles are reviewed in section \ref{Poles}).

\vskip 0.2 cm 
- Taylor series.

\vskip 0.2 cm 
- Logarithmic behaviours, namely series expansions  with coefficients which are polynomials of $\ln x$. $y(x)\to$ constant as $x\to 0$, where the constant may be zero or not. 

\vskip 0.2 cm 
- Inverse  logarithmic behaviours,   expanded as formal series of $( \ln x)^{-1}$. $y(x)=O(1/\ln x)$, or $O(1/\ln^2 x)$,  as $x\to 0$. 

\vskip 0.2 cm 
For any $\vartheta>0$ there exists  $r(\vartheta)$  small enough (decreasing function of $\vartheta$) such that the expansion of the complex power behaviours converges for  $|\arg(x)|<\vartheta$ and $0<|x|<r(\vartheta)$. If  $\Re\sigma=0$, there is also an  additional constraint $\arg (x)>\varphi_0$, for a suitable $\varphi_0$ fixed by the integration constants.  Also the expansion of the denominator of the inverse oscillatory behaviours is convergent under  the same conditions $|\arg(x)|<\vartheta$, for any choice of a  $\vartheta>0$, and $0<|x|<r(\vartheta)$, plus the additional contraint $\arg (x)<\varphi_0$, for a suitable $\varphi_0$ (see section \ref{Poles} and \ref{conver}).     When logarithms appear, no proof of convergence is known yet. Convergence of the Taylor expansions is studied in \cite{Kaneko}. 

 In the table, $\sigma,\phi,\nu$ and $a$ denote integration constants. The coefficients $c_{nm}$, $d_{nm}$ and  $b_n$ are 
rational functions of $\sqrt{\alpha}$, $\sqrt{\beta}$, $\sqrt{\gamma}$ and $\sqrt{1-2\delta}$. The coefficients $b_n(a)$  are 
rational functions of $\sqrt{\alpha}$, $\sqrt{\beta}$, $\sqrt{\gamma}$, $\sqrt{1-2\delta}$ and $a$. $P_n(\ln x, a)$ are 
polynomials in $\ln x$ with coefficients which are rational functions of $\sqrt{\alpha}$, $\sqrt{\beta}$, $\sqrt{\gamma}$, 
$\sqrt{1-2\delta}$ and $a$. These coefficients can be recursively computed (essentially, by substitution into PVI).

$$
\begin{array}{|| c|c|c||} 
\hline\hline   &  &  \cr 
  \hbox{{\bf  Complex power behaviours}} & \matrix{\hbox{Integ.}\cr
                                                   \hbox{ const.}
                                                   }
 & \hbox{Other Conditions}    \cr
\hline          
 & &  \cr
\matrix{
y(x)= \sum_{n=1}^\infty x^n\sum_{m=-n}^n c_{nm}(ax^{\sigma})^m 
\cr
\cr
\hbox{
This is (\ref{locintro}) and 
(\ref{hoshi})}. }
& \matrix{\sigma \cr
\cr  a\neq 0} & 
 \matrix{
0\leq \Re \sigma<1,
\cr
2\cos\pi\sigma=p_{0x},
\cr
p_{0x}\not\in(-\infty,-2]\cup\{2\}
}
 \cr
&&\cr
\hline
&&\cr
%
%
y(x)={\sqrt{\alpha} +(-)^k\sqrt{\gamma} \over \sqrt{\alpha}}+ \sum_{n=1}^\infty x^n \sum_{m=0}^n d_{nm}(\tilde{a}x^{\rho})^m
&
a
&
\matrix{
\alpha\neq 0.
\cr
p_{0x}=2\cos\pi\rho\neq \pm 2.
\cr 
\rho+1=
\cr
=(\sqrt{2\alpha}\pm\sqrt{2\gamma})\hbox{sgn}(\sqrt{2\alpha}\pm\sqrt{2\gamma})
}
\cr
 && \cr 
\hline && \cr 
%
%
y(x)={ 1\over a }x^{-\omega}
\left(
1+\sum_{n=1}^\infty x^n \sum_{m=0}^n d_{nm}(a x^{\omega})^m
\right)
&
a
&
\matrix{\alpha =0.
\cr 
\omega=\sqrt{2\gamma}~\hbox{sgn}(\Re\sqrt{2\gamma}),
\cr
\Re\omega > 0.
}
\cr 
&&
\cr
\hline
\end{array}
$$


$$
\begin{array}{|| c|c|c||} 
\hline\hline
&&\cr
\hbox{{\bf Inverse Oscillatory Behaviours}} & \left.\matrix{\hbox{Integration}\cr\hbox{constants}}\right. & \hbox{ Other }
\cr
\hline
&&
\cr
\matrix{
y(x)= \Bigl[\sum_{n=0}^\infty x^n\sum_{m=-n-1}^{n+1} c_{nm}
\bigl(e^{i\phi}x^{2i\nu}\bigr)^m
\Bigr]^{-1} 
\cr
\cr 
= \left[A\sin(2\nu \ln x +\phi)+B+ O(x)\right]^{-1}
\cr 
\cr
A=-\sqrt{{\alpha\over 2\nu^2}+B^2},~~
~~~B={2\nu^2+\gamma-\alpha\over 4\nu^2}
}
                          &
 \matrix{\nu 
\cr 
\cr 
\phi}
                             &
\matrix{
\nu\in{\bf R}\backslash\{0\},
\cr 
-2\cosh2\pi\nu=p_{0x}<-2.
}
\cr
&&
\cr 
\hline
&&\cr
%
%
y(x)=\left[\sum_{n=0}^\infty x^n
\sum_{m=0}^{n+1}c_{n+1,m}\bigl(ax^{-2i\nu}\bigr)^m
\right]^{-1}
&
a
&
\left.
\matrix{2i\nu= \pm(\sqrt{2\alpha}\pm\sqrt{2\gamma})
\cr
\in i {\bf R}\backslash\{0\}.
}
\right.
\cr
&&\cr
\hline
\end{array}
$$



$$
\begin{array}{||c|c|c||}
\hline\hline
&&\cr
 \hbox{{\bf Taylor expansions}} & \left.\matrix{\hbox{Int.}\cr\hbox{const.}}\right. & \hbox{ Other Conditions}
\cr
\hline
&&
\cr
 y(x)={\sqrt{-2\beta}\over \sqrt{-2\beta}+(-)^k\sqrt{1-2\delta}}x+\sum_{n=2}^\infty b_nx^n   
& &
                  \matrix{ \sqrt{-2\beta}\pm\sqrt{1-2\delta}\not\in{\bf Z}.
\cr 
\cr
 p_{0x} \neq \pm 2.
}
\cr
&& \cr
\hline
%
%
&&
\cr
y(x)= \sum_{n=1}^{|N|}b_nx^n+ax^{|N|+1}+\sum_{n= |N|+2}^\infty b_n(a)x^n
&
a
&
                                        \matrix{\sqrt{-2\beta}\pm\sqrt{1-2\delta} =N\neq 0
\cr
  p_{0x}=2\cos \pi N=\pm2.
}
\cr
&&\cr
\hline
&&\cr      
                                                       y(x)= ax+a(a-1)\left(\gamma-\alpha-{1\over 2}\right)x^2+\sum_{n= 3}^\infty b_n(a)x^n
&
        a
& 
\matrix{2\beta=2\delta-1=0.
\cr 
 p_{0x}=2.
}
\cr
&&\cr
\hline
&&\cr
%
%
y(x)= {\sqrt{\alpha} +(-)^k \sqrt{\gamma} \over \sqrt{\alpha}} 
+
\sum_{n=1}^\infty b_nx^n,~~\hbox{ Basic Taylor}
&
&
   \matrix{                            \alpha\neq 0, 
\cr
            \sqrt{2\alpha}+(-)^k\sqrt{2\gamma}\not\in{\bf Z}.
\cr 
 p_{0x}\neq \pm 2.
}
  \cr
&&
\cr
\hline
&&\cr
%
%
y(x)= \sum_{n=0}^{|N|-1}b_n x^n+ax^{|N|}+\sum_{n= |N|+1}^\infty b_n(a)x^n
&
    a  &
            \matrix{\sqrt{2\alpha}\pm\sqrt{2\gamma} =N \neq 0
\cr
p_{0x}=-2\cos\pi N=\pm 2.
}
\cr
&&\cr
\hline
&&\cr
y(x)=a+(1-a)(\delta-\beta)x+\sum_{n= 2}^\infty b_n(a)x^n
      &
           a
                 &
\matrix{\alpha=\gamma=0.
\cr 
 p_{0x}=-2.
}
\cr
&&\cr
\hline
\end{array}
$$


$$
\begin{array}{|| c|c|c||}
\hline\hline
&&\cr
\hbox{{\bf Logarithmic behaviours}}
 & \left.\matrix{\hbox{Int.}\cr\hbox{const.}}\right. & \hbox{Other Conditions}
\cr
\hline
&&
\cr
\matrix{
y(x)=\Sigma_{n=1}^{|N|}b_nx^n
+\Bigl(a+b_{|N|+1}\ln x\Bigr)x^{|N|+1}+ 
\cr 
\cr 
+\Sigma_{n= |N|+2}^\infty P_n(\ln x;a)x^n,~~N\neq 0.
\cr
\cr
 y(x)= \Bigl(a\pm\sqrt{-2\beta}\ln x\Bigr)x+\sum_{n=2}^\infty P_n(\ln x;a)x^n,~~N= 0.
} 
                                                               & 
a
                                                               &
       \matrix{
\sqrt{-2\beta}\pm\sqrt{1-2\delta} =N
 \cr
p_{0x}=2\cos\pi N=\pm 2
}
\cr
&&\cr
\hline
&&\cr
\matrix{y(x)=\left[{2\beta+1-2\delta\over 4}(a+\ln x)^2+{2\beta\over 2\beta+1-2\delta}\right]x+
\cr
\cr
+\sum_{n\geq 2}^\infty P_n(\ln x;a)x^n }
& a & 
\matrix{2\beta\neq 2\delta-1.
\cr
p_{0x}=2.
}
\cr
&&\cr
\hline
&&\cr
        \matrix{
  y(x)= \sum_{n=0}^{|N|-1}b_n x^n+\Bigl(a+b_N\ln x\Bigr)x^{|N|}+
\cr 
\cr 
+\sum_{n= |N|+1}^\infty P_n(\ln x;a)x^n
}
& 
a
&
                                \matrix{
                        \sqrt{2\alpha}\pm\sqrt{2\gamma} =N\neq 0
\cr
p_{0x}=-2\cos \pi N=\pm 2.
}
\cr
&&
\cr
\hline
\end{array}
$$


$$
\begin{array}{|| c|c|c||} 
\hline\hline
&&\cr
\hbox{{\bf Inverse logarithmic behaviours}}
& \matrix{\hbox{Int.}
\cr 
\hbox{const.}
} & \hbox{ Other Conditions }
\cr
\hline
&&
\cr
\left.
\matrix{
y(x)= \left\{a\pm \sqrt{2\alpha}\ln x+\sum_{n=1}^\infty P_n(\ln x;a)x^n\right\}^{-1}
\cr
\cr
= \pm{1\over \sqrt{2\alpha}\ln x}\left[1\mp {a\over \sqrt{2\alpha} \ln x}+O\left({1\over \ln^2 x}\right)\right]
}
\right.
&
a
&\matrix{\alpha=\gamma\neq 0.
\cr 
p_{0x}=-2.
}
\cr
&&\cr
\hline
&&\cr
\left.
\matrix{
y(x)= \left\{{\alpha\over \alpha -\gamma}+ {\gamma-\alpha\over 2} (a+\ln x)^2 +\sum_{n=1}^\infty P_{n+1}(\ln x;a) x^n\right\}^{-1}
\cr\cr
=
{2\over (\gamma-\alpha)\ln^2 x}\left[
1-{2a\over \ln x} +O\left({1\over \ln^2 x}\right)
\right]
}
\right.
&
a
&
\matrix{
\alpha\neq \gamma.
\cr 
p_{0x}=-2 .
}
\cr
&&\cr
\hline
\end{array}
$$


\subsection{The Parametric Connection Formulae (\ref{POA}) and (\ref{POA1})}
 The parametric connection formulae are quite long to write, so  we refer  to \cite{guz2012}, which collects  the 
 formulae computed in  \cite{Jimbo}, \cite{Boalch},  \cite{D4}, \cite{D3},  \cite{D2}, \cite{D1}, 
\cite{guz2010}. 
  Just to give the taste of them, we give a very simple example. Consider the following solution with Taylor expansion at $x=0$ (namely, $u=0$): 
$$y(x)=ax+a(a-1)\left(
 \gamma-\alpha-{1\over 2}
\right)x^2+O(x^3).
$$
 The formulae (\ref{POA}) are 
$$
a = 
 {2\cos\pi\sqrt{2\gamma}-p_{01}
\over 
4\cos\left({\pi\over 2}\left[\sqrt{2\alpha}+\sqrt{2\gamma}  \right] \right)
     \cos\left({\pi\over 2}\left[\sqrt{2\alpha}-\sqrt{2\gamma}  \right]\right)
}
$$
and the (\ref{POA1}) are
$$
p_{0x}=2,
$$
$$ 
p_{01}=2\cos\pi\sqrt{2\gamma}-4a\cos\left({\pi\over 2}\left[\sqrt{2\alpha}+\sqrt{2\gamma}  \right] \right)
     \cos\left({\pi\over 2}\left[\sqrt{2\alpha}-\sqrt{2\gamma}  \right]\right),
$$ 
$$ 
p_{x1}=2\cos\pi\sqrt{2\gamma}+4(a-1)\cos\left({\pi\over 2}\left[\sqrt{2\alpha}+\sqrt{2\gamma}  \right] \right)
     \cos\left({\pi\over 2}\left[\sqrt{2\alpha}-\sqrt{2\gamma}  \right]\right).
$$
In the above, $p_\mu$ (or $\theta_\mu$), $\mu=0,x,1,\infty$ are re-expressed in terms of the coefficients of $\hbox{PVI}_{\alpha\beta\gamma\delta}$.

\subsection{Other Approaches to obtain Critical Behaviours} 
Two alternative approaches  to the local analysis of PVI provide the critical behaviours, with full expansions. 

The first is a local analysis of integral equations associated to PVI, due to Shimomura \cite{Sh}, \cite{IKSY}. It provides  behaviours of type (\ref{locintro}) and (\ref{hoshi0}),  on the universal covering of a punctured neighbourhood of the critical point. The elliptic representation of PVI gives the same results \cite{D4} \cite{D3}. We will come back to this point in Section \ref{conver}. 

The second approach is the method of {\it power geometry}  developed by Bruno  in \cite{Bruno4} and  \cite{Bruno5}.  This method is a refinement of the method of Newton Polygons. In  a series of papers  summarised in the review \cite{Bruno7}, Bruno and  Goryuchkina construct all the (possibly formal) expansions that can be obtained by means of this method. They are classified in a way which resents of the method.  In \cite{guz2012}  we  proved  that they coincide with those in the table reviewed  in section \ref{tabulation}.

We remark that the local analysis does not provide the connection formulae, for which we need to go back to the isomonodromy deformation method.

\section{Poles}
\label{Poles}
The position of the poles is known for classical functions, but  for $\hbox{PVI}_{\alpha\beta\gamma\delta}$ we still  do not know {\it the global  distribution of the poles}, except for  special choices of $\alpha,\beta,\gamma,\delta$. For example, in \cite{Br}, the pole distribution for  $\hbox{PVI}_{{1\over 8}{1\over 8}{1\over 8}{3\over 8}}$ (the Hitchin equation \cite{hitchin}) is determined  on the whole universal covering of
 ${\bf C}\backslash\{0,1,\infty\}$. A formula for an infinite series of poles is given in terms of Theta-functions. The poles are distributed along  lines which  are spirals at a small scale  around the critical points, and more complicated lines on the whole universal covering. A birational Okamoto's transformation  transforms $\hbox{PVI}_{{1\over 8}{1\over 8}{1\over 8}{3\over 8}}$  into $\hbox{PVI}_{000{1\over 2}}$, known as  Picard equation \cite{Picard}.

 Though the { global distribution of the poles} is limited to the above examples, we know  the {\it asymptotic distribution of the poles close to a critical point}.  The existence of poles close to the critical point $x=0$ is due the the structure of the branches  (\ref{invos1}), because the denominator may vanish when $x$ approaches zero. The same may happen at $x=1$ and $x=\infty$ for (\ref{invos2}) and (\ref{invos3}) respectively. To understand the  asymptotic distribution of these  poles, let us concentrate on the critical point $x=0$, being the others analogous.  

The branch  (\ref{invos1}) is obtained through (\ref{sym2}) from the branch (\ref{hoshi0}) = (\ref{hoshi}). As discussed in section \ref{JimboWork}, (\ref{locintro}) and  (\ref{hoshi0})  hold when $0<|x|<r$ and $|\arg x|<\vartheta$, where $\vartheta>0$ is a positive number (it can be chosen arbitrarily), and    $r=r(\vartheta)$ is  small and decreasing  when $\vartheta$  increases.  The symmetry (\ref{sym2}) does not affect the condition  $|\arg x|<\vartheta$, $0<|x|<r$,  therefore  (\ref{invos1}) holds under the same condition as well.  
Since $\arg x$ remains bounded,  $x$ may tend to zero along {\it a radial path, while spiral paths are not allowed}. With this restriction, most solutions of PVI have no poles in a sufficiently small neighbourhood of $x=0$, except precisely  the class of solutions   (\ref{invos1}). Let us denote a solution  (\ref{invos1}) with  $y=y(x,\nu,\phi)$, and re-write:  
$$
 y(x,\nu,\phi)= {1\over y_1(x)+g(x)},
$$
$$
y_1(x):= {\cal A}\sin(2\nu~\ln x +\phi)+{\cal B},~~~~~g(x)=O(x),
$$
It follows that $y(x,\nu,\phi)$   has poles in a neighbourhood of $x=0$, coinciding with the zeros of $y_1(x)+g(x)$. They are asymptotically close to the zeros of $y_1(x)$ when $x\to 0$, because $y(x,\nu,\phi)^{-1}\sim y_1(x)$.  Direct computation shows that  $y_1(x)$ has  two infinite sequences of zeros,  distributed along two rays converging to $x=0$ (see figure \ref{POLI}).

The full expansion of (\ref{invos1}) is in the table of section \ref{tabulation}. It is convenient to write it in reciprocal form:
\be
\label{pro2y}
{1\over y(x,\nu,\phi)} = y_1(x)+xy_2(x)+x^2y_3(x)+...~=\sum_{n=1}^\infty x^{n-1} y_n(x)
\ee
where: 
\be
\label{coefA}
 y_n(x)=\sum_{m=-n}^n A_{nm}(\nu,\alpha,\beta,\gamma,\delta)~ e^{im\phi}x^{2im\nu},~~~0\leq \Re\phi\leq \pi. 
\ee
The $A_{nm}(\nu, \alpha,\beta,\gamma,\delta)$'s are  algebraic functions of $\nu,\alpha,\beta,\gamma,\delta$. 
 Their explicit form is recursively  computable  by the procedure of \cite{guz2010}.
 The series (\ref{pro2y}) can be formally computed for any $x$. There exists  $r<1$ small enough such that it   converges 
(at least) in the domain 
\be\label{appe}
\ln|x|-\Im\phi-\ln r <2\nu \arg(x)<-\Im \phi+\ln r,~~~~~0<|x|<r.
\ee
This result is reviewed in section \ref{conver}.\footnote{
Note that, in general, for a series like (\ref{pro2y}), one expects convergence for $0<|x|<r$ and $|x^{n+2im\nu}e^{im\phi}|<\epsilon_{nm}$, where $r,~\epsilon_{nm}>0$ are sufficiently small. Thus
$$
\ln|x|-\Im\phi +\max\Bigl\{-\ln r,\sup_{m>0,n\geq 1}\left|
{\ln \epsilon_{nm}\over m}
\right|\Bigr\}
<
2\nu\arg x
<
-\Im\phi+\min\Bigl\{\ln r,\inf_{m<0,n\geq 1}\left|
{\ln \epsilon_{nm}\over m}
\right|\Bigr\}.
$$
}  
The $n$-th order  $y_n(x)$ is an oscillatory bounded function for $x\to 0$ and $\arg x$ bounded, and  $x^{n-1}y_n(x)=O(x^{n-1})$. We remark that (\ref{pro2y}) does not vanish in the above domain. Thus, the poles of $y(x,\nu,\phi)$ necessarily lie in the domain $2\nu \arg x\geq -\Im \phi+\ln r$. 
The following theorem estabilishes the final result \cite{Dpoli}: 
\begin{figure}
\centerline{\includegraphics[width=0.6\textwidth]{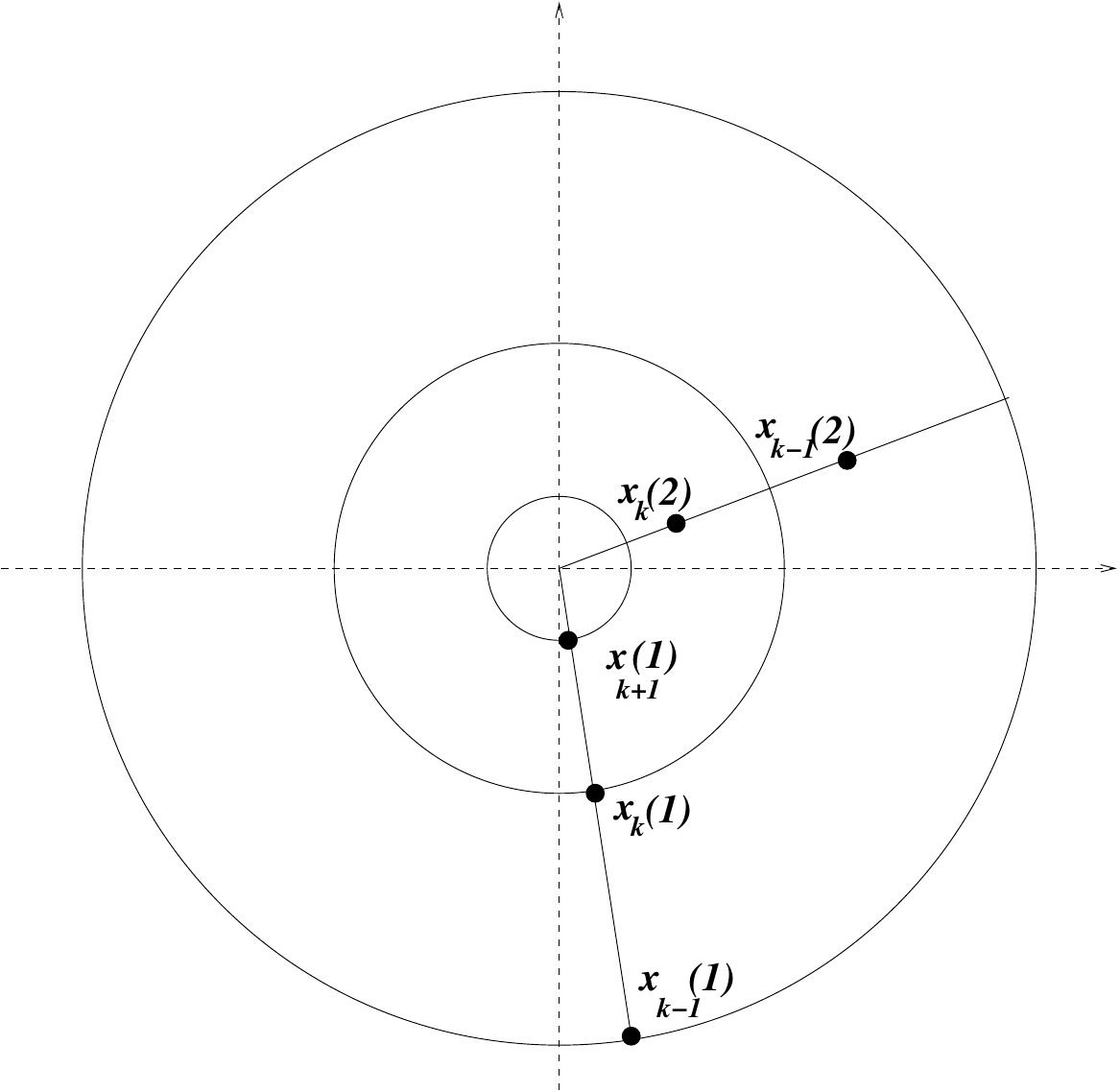}}
\caption{The zeros of $y_1(x)$ around $x=0$.}
\label{POLI}
\end{figure}
\bth
\label{theo1}
 Let 
$\nu>0$ and let $y(x,\nu,\phi)$ be (\ref{invos1}).  
  There are two sequences $\{x_k(1)\}_{k\in{\bf Z}}$ and $\{x_k(2)\}_{k\in{\bf Z}}$ of zeros of $y_1(x)$:
\be
\label{seRO}
x_k(j)=\exp\left\{
-{\phi\over 2\nu} -{i\over 2\nu} \ln \left[
(-)^j \sqrt{{A_{10}^2\over 4 A_{11}^2}-{A_{1,-1}\over A_{11}}}-{A_{10}\over 2 A_{11}} 
\right]
-{k\pi\over \nu}\right\}
,~~~k\in{\bf Z},~~~j=1,2 
\ee
The argument of $\ln[~]$ is fixed, being other choices  absorbed into  ${k \pi \over \nu}$.
Let $k_0\in{\bf N}$  be sufficiently big in order to have $|x_k(j)|<r(\vartheta)$, $j=1,2$, in such a way that (\ref{invos1}) holds.  There exists $K$ sufficiently big such that for every  $k\geq \max\{K,k_0\}$, and every $j=1,2$, $y(x,\nu,\phi)$ has a  pole $\xi_k(j)$ lying in a neighbourhood of $x_k(j)$,  with   the  asymptotic representation
\be
\label{POLESintro}
   \xi_k(j)=x_k(j)-{1\over 2} x_k(j)^2 +\sum_{N=3}^\infty\Delta_N(j) x_k(j)^N, ~~~k\to +\infty,~~~x_k(j)\to 0.
\ee
The coefficients $\Delta_N(j)\in{\bf C}$ are certain numbers independent of $k$  that can be computed form the coefficients $A_{nm}$ of (\ref{coefA}).
 There are no other poles  for $|x|$ sufficiently small. 
\eth

The series (\ref{POLESintro}), computed making use of (\ref{pro2y}) and (\ref{coefA}), is at least asymptotic. We can actually prove its convergence in some cases, for example (see \cite{Dpoli}) when the equation is $\hbox{PVI}_{000{1\over 2}}$,  the already mentioned Picard equation.
The position of the poles depends on 
 $\nu$ and $\phi$, therefore it depends on the monodromy data.

\vskip 0.2 cm
\noindent
{\bf Remark:} Let $\vartheta=\pi$. The zeros $x_k(j)$ and the poles may be outside the range $|\arg x|<\pi$, depending on $\Im \phi$. However, 
 the analytic continuation of the branch $y(x,\nu,\phi)$ corresponding to the loop $x\mapsto xe^{2\pi i}$ is \cite{Dpoli}:
\be
\label{analcont}
y(x,\nu,\phi)\mapsto   y(x,\nu,\phi+2\pi i\nu)=y(xe^{2\pi i},\nu,\phi),~~~-\pi\leq \arg x <\pi.
\ee
The shift $\phi  \mapsto \phi+2\pi i\nu$ changes the imaginary exponent of the $x_k(1)$'s by $-2\pi i$. This implies that, by a sufficient number of loops,  we can always find a branch with    at least one  of the two sequences of poles   in the range $-\pi \leq \arg x<\pi$.

 \vskip 0.2 cm

\noindent
{\bf Remark:} (\ref{seRO}) implies that 
$$
2\nu\arg x_k(j)=-\Im \phi -\ln \left|
(-)^j \sqrt{{A_{10}^2\over 4 A_{11}^2}-{A_{1,-1}\over A_{11}}}-{A_{10}\over 2 A_{11}} 
\right|.
$$
Moreover, we know that the poles exist if $2\nu \arg x\geq -\Im \phi+\ln r$. Thus, we obtain an upper estimate of the radius of convergence of  (\ref{pro2y}):
$$
 r\leq \hbox{min}\left\{e^{
-
\ln \left|
\sqrt{{A_{10}^2\over 4 A_{11}^2}-{A_{1,-1}\over A_{11}}}-{A_{10}\over 2 A_{11}} 
\right|
},~e^{
-
\ln \left|
 \sqrt{{A_{10}^2\over 4 A_{11}^2}-{A_{1,-1}\over A_{11}}}+{A_{10}\over 2 A_{11}} 
\right|
}
\right\}.
$$

 \vskip 0.2 cm

\noindent
{\bf Example:} If $\alpha={9\over 2}$,   $\beta=\gamma=0$ and $\delta={1\over 2}$, then (\ref{pro2y}) is  a  branch associated to the {\it quantum cohomology of $ CP^2$} \cite{Dub1}, \cite{Dub2}.  As it is proved in \cite{Dpoli}, the branch is identified by the following integration constants
$$
\nu = {2 \ln {\bf G}\over \pi}=0.30634...
$$
where ${\bf G}={1+\sqrt{5}\over 2}$ is the golden ratio (I thank M. Mazzocco for this remark), and 
\be
\label{3vetrini}
\phi=i \ln
\left\{-
\pi^2~{({\bf G}^4+1)^2\over ({\bf G}^2+1)^2}~2^{16i\nu}~{(1-2i\nu)^2 ~\nu^2 \over (1+2i\nu)^2}~{\Gamma(1-2i\nu)^4\over 
\Gamma(1-i\nu)^8}
\right\},
\ee
The first approximation $y_1(x)$ has  two  infinite sequences of zeros  accumulating at $x=0$ along {\it the negative imaginary axis}, because ${\Im \phi\over\nu}={\pi\over 2}$:   
\be
\label{synzeros}
  x_k(j)= -i\exp
\left\{
-{\Re \phi\over 2\nu} -{2(j-1)\over \nu}\left|\arccos{3\over \sqrt{4\nu^2+9}}\right|\right\}\exp\left\{-{k\pi\over \nu}\right\}, ~~k\in{\bf N}, ~~~j=1,2.
\ee
The first coefficients in the asymptotic expansion (\ref{POLESintro}) are 
$$
\Delta_3(1)={176\nu^4+185+352\nu^2 \over 1024 (\nu^2+1)^2}=0.1792...,~~~
\Delta_3(2)={401+176\nu^4+352\nu^2\over 1024(\nu^2+1)^2}=0.3555...
$$
$$
\Delta_4(1)=-{57+48\nu^4+96\nu^2\over 1024(\nu^2+1)^2}=-0.05422...,~~~
\Delta_4(2)=-{273+48\nu^4+96\nu^2)\over 1024(\nu^2+1)^2}=-0.2305...
$$

\section{Critical Behaviours on  the Universal Covering of a Punctured Neighborhood of a Critical Point}
\label{conver}
 Let us consider again the critical point $x=0$. 
 The expansions of the complex power behaviours and the inverse oscillatory behaviours,  in the table of section \ref{tabulation}, converge in suitable domains. This is a consequence of the fact that such expansions (basically the expansions of  (\ref{locintro}), (\ref{hoshi0}) and (\ref{invos1})), coincide with the convergent expansions obtained by a different local analysis of PVI, which we are going to discuss in subsections \ref{shimomura} and \ref{eliptic}. 

\subsection{Shimomura's Representation }
\label{shimomura}
  The approach to the local analysis of critical behaviours, due to S. Shimomura \cite{Sh} \cite{IKSY}, consists in rewriting PVI as a system of two differential equations of the first order, and then as a system of two integral equations. The latter is solved by successive approximations and the solutions are obtained  in the form of  convergent expansions, as stated in the following theorem:

\bth
\it Let $\widetilde{{\bf C}_0}$ be the universal covering of a punctured neighbourhood of $x=0$.  For any complex number $a\neq 0$ and for any $ \sigma \not \in   (-\infty,0]
\cup[1,+\infty)$  there is a 
sufficiently small $r$  such that $\hbox{PVI}_ {\alpha\beta\gamma\delta}$  has a 
holomorphic solution in the domain 
$$
  {\cal D}_s(r;\sigma,a)= \{{x} \in \widetilde{{\bf C}_0} ~|~ |x|<r, 
~|a{x}^{1-\sigma}|<4r,~|a^{-1} {x}^{\sigma}|<r/4 \}
$$
with the following representation:
\be
\label{shimo}
   y(x,\sigma,a)= {1 \over \cosh^2({\sigma-1\over 2}\ln {x}
    -{1\over 2}\ln a+\ln 2 +{v({x},\sigma,a)\over 2})},
\ee
where 
\be 
\label{typeseries}
  v({x},\sigma,a)= \sum_{n\geq 1} a_n(\sigma) {x}^n+ \sum_{n\geq 0, 
~m\geq 1} b_{nm}(\sigma) {x}^n (a{x}^{1-\sigma})^m 
+ \sum_{n\geq 0,~m\geq 1} c_{nm}(\sigma) {x}^n 
(a^{-1}{x}^{\sigma})^m,
\ee
is convergent 
(and holomorphic) in ${\cal D}_s(r;\sigma,a)$. Here 
 $a_n(\sigma), b_{nm}(\sigma), c_{nm}(\sigma)$ are certain  rational functions of 
$\sigma$.  
\eth

\noindent
The denominator of (\ref{shimo}) does not vanish in $ {\cal D}_s(r;\sigma,a)$. The domain ${\cal D}_s(r;\sigma,a)$ is an open domain in the plane $(\ln|x|,\arg(x))$:
\be
 |x|<r,~~~ \Re\sigma \ln|x| -\ln |a|  - \ln {r\over 4} <
                                                    \Im \sigma \arg(x)
           < (\Re \sigma -1 )\ln |x| -\ln |a|  + \ln (4r).
\label{UFFAUFFA}
\ee
  If $\Im \sigma=0$, the domain is simply ${\cal D}_0(r):= 
\left\{ x\in \widetilde{{\bf C}_0}~ \hbox{ such that }
  |x|<r \right\}$.  
We can  compute the critical behaviour of (\ref{shimo})  for $x\to 0$ along a regular path contained in  $ {\cal D}_s(r;\sigma,a)$, connecting a point $x_0\in {\cal D}_s(r;\sigma,a)$ to $x=0$. Let us consider the   one parameter ($\Sigma$) family of paths
\be
\label{PATH}
|x|<r,~~~\arg x= \arg x_0 +{\Re \sigma-\Sigma\over \Im \sigma}(\ln|x|-\ln |x_0|),~~~~~0\leq \Sigma\leq 1
\ee
or any regular path if $\Im \sigma=0$.  If  $0\leq \Sigma< 1$, then (\ref{shimo}) behaves as follows
\be
\label{locintro1}
 y(x,\sigma,a)=ax^{1-\sigma}(1+O(x+x^{1-\sigma}+x^\sigma)), ~~~~~\hbox{ if } ~~0<\Sigma<1,
\ee
or
\be
\label{hoshi1}
y(x,\sigma,a)= x\sin^2\left({i\sigma\over 2}\ln x  + \psi(x)\right)+O(x^2),~~~~~\hbox{ if }~~ \Sigma=0.
\ee
The function $\psi(x)$ has a convergent expansion (following from $v(x)$) of type 
$$
\psi(x)=\sum_{n\geq 0} \psi_n (\sigma, a)x^{n\sigma},
$$
 for certain $\psi_n(\sigma,a)$'s. It oscillates without vanishing when $x\to 0$ along the path with  $\Sigma=0$. 
It is clear that (\ref{locintro1}) coincides with (\ref{locintro}) when $0\leq \Re\sigma<1$ and $x\to 0$ along a radial path, with the identification  $\sigma_{0x}=\sigma$.    On the other hand, the  behaviour (\ref{hoshi1}) can be rewitten 
as 
\be
\label{HOshi1}
y(x,\sigma,a)= x\Bigl\{A\sin(i\sigma \ln x+\phi)+B\Bigr\}+O(x^2),~~~~~\Sigma=0.
\ee
with the same coefficients $A$ and $B$ of (\ref{hoshi0}). This follows from the fact that we can write (see \cite{guz2010})
$$
A\sin(i\sigma \ln x+\phi)+B=-2A\sin^2\left({i\sigma\over 2}\ln x+{\phi\over 2}-{\pi\over 4}\right)+A+B=
$$
$$
= \sin^2\left({i\sigma\over 2}\ln x +\sum_{n\geq 0} \psi_nx^{n\sigma}\right) ~~~\hbox{ for suitable }\psi_n.
$$
Thus, (\ref{hoshi1}) 
coincides with (\ref{hoshi0}) when $\Re\sigma=0$ and $x\to 0$ along a radial path. 
When $\Sigma=1$ we obtain the critical behaviour
\be
\label{invasi}
y(x,\sigma,a)={1\over \sin^2\left(i{\sigma-1\over 2}\ln x  + \psi(x)\right)+O(x)},~~~~~\Sigma=1.
\ee
The function $\psi(x)$ has a convergent expansion 
$$
\psi(x)=\sum_{n\geq 0} \psi_n (\sigma, a)x^{n(1-\sigma)},
 $$ 
for certain $\psi_n(\sigma,a)$, which does not vanish along the path with $\Sigma=1$.   
By the same argument above, $y(x,\sigma,a)$ coincides with
\footnote{~~~~$ 
{\cal A}^2={\bigl[(1-\sigma)^2-(\theta_\infty-1-\theta_1)^2\bigr]\bigl[(1-\sigma)^2-(\theta_\infty-1+\theta_1)^2\bigr]
\over 4(1-\sigma)^2},
$ ~~~
$
{\cal B}={(\theta_\infty-1)^2-\theta_1^2+(1-\sigma)^2\over 2(1-\sigma)^2}
$. }
\be
\label{invos11}
y(x)={1\over {\cal A}\sin(i(1-\sigma)\ln x+\phi)+{\cal B}+O(x)}
, ~~~~~x\to 0.
\ee
 If $\Re\sigma=1$, then (\ref{invos11})  {\it coincides} with (\ref{invos1}) (with $\sigma=1+2i\nu$, $\nu\in{\bf R}$), and the convergence $x\to 0$ is along a radial path with $\arg x =\arg x_0$. 

 It is crucial  to observe  that from (\ref{shimo}) and (\ref{typeseries}) we can compute the {\it full expansions} {\it convergent on the domains} ${\cal D}_s(r;\sigma,a)$,  of (\ref{locintro1}), (\ref{hoshi1}) and of the denominator of (\ref{invasi}).  The one to one correspondence between a branch of a transcendent and a point of the space of monodromy data - when holds -  implies that only one transcendent has the given critical behaviour.  
As a result of this, the  convergence of the full expansion corresponding to the behaviours (\ref{locintro}), (\ref{hoshi0}) and (\ref{invos1}) is proved.

  Shimomura's result applies to any value of $\sigma_{0x}$, except for values in $(-\infty,0]
\cup[1,+\infty)$, provided that the convergence $x\to 0$ occurs along spirals whenever $\Re\sigma_{0x}<0$ or $\Re\sigma_{0x}>0$. Therefore, it is  a generalization of  (\ref{locintro}), (\ref{hoshi0})  and (\ref{invos1}), which hold only for $0\leq \Re\sigma_{0x}\leq 1$. 

 Finally, we recall that the above results have been used in section \ref{Poles}, where the domain (\ref{appe}) is the particular case of ${\cal D}_s(r;\sigma,a)$ when   $\sigma=1+2i\nu$, $\nu\in{\bf R}$, that the denominator of (\ref{invos1})= (\ref{invos11}) does not vanish in (\ref{appe}), and the poles are outside this domain.
\vskip 0.2 cm 


\subsection{The Elliptic Representation}
\label{eliptic}

The elliptic representation was introduced by  P. Painlev\'e in 
 \cite{PaiCRAS1906} and R. Fuchs in \cite{fuchs}. Let $\wp(z;\omega_1,\omega_2)$ be 
 the Weierstrass Elliptic function of the independent variable $z\in 
{\bf P}^1$, with 
 {\it half-periods } $\omega_1$, $\omega_2$. Let us consider the following $x$ dependent half periods
$$
\omega_1(x):={\pi \over 2} F\left({1\over 2},{1\over 2},1;x\right), ~~~~
\omega_2(x):=i{ \pi \over 2}F\left({1\over 2},{1\over 2},1;1-x\right), 
$$
where $F$ are hypergeometric functions. 
A PVI transcendent can be represented as follows \cite{fuchs}:  
$$
   y(x)= \wp \left({u(x)\over 2};\omega_1(x),\omega_2(x) \right)
+{1+x\over 3}. 
$$ 
where $u(x)$ solves a non linear  ODE of the 2nd order, equivalent to PVI.
  The algebraic-geometrical properties of the elliptic representations  
where studied   
in \cite{Manin2}. The  analytic properties of the function $u(x)$ in the simplest case 
  $\alpha=\beta = \gamma=1-2\delta=0$ are easily studied, because the equation for $u(x)$ becomes a Gauss hypergeometric equation. Accordingly, the function $u(x)$ is a 
linear combination of $\omega_1$ and $\omega_2$. This case was well known 
to Picard  \cite{Picard}, and we have already mentioned it. 

The local (at $x=0$) analytic properties of $u(x)$  were studied in the general case in \cite{D4} and \cite{D3}, as follows.  
Let $\nu_1$ and $\nu_2$ be complex numbers. 
 Let us  consider the domains 
 \be
{\cal D}(r;\nu_1,\nu_2):= 
\left\{ x\in \widetilde{{\bf C}_0}~ \hbox{ such that }
  |x|<r, \left|{e^{-i\pi \nu_1}\over 16^{1-\nu_2}} x^{1-\nu_2} 
\right|<r,
\left| 
{e^{i\pi \nu_1} \over 16^{\nu_2}} x^{\nu_2}\right|<r \right\},
\label{DOMINOO1talk}
\ee
\be
{\cal D}_0(r):= 
\left\{ x\in \widetilde{{\bf C}_0}~ \hbox{ such that }
  |x|<r \right\},
\label{DOMINOO3}
\ee 
and the expansion:
\be 
v(x,\nu_1,\nu_2):= \sum_{n\geq 1} a_n x^n + \sum_{n\geq 0, m\geq 1} b_{nm} x^n 
\left[ e^{-i\pi \nu_1} \left({x\over 16}\right)^{1-\nu_2}\right]^m 
+\sum_{n\geq 0,m\geq 1} c_{nm} x^n 
\left[ e^{i\pi \nu_1} \left({x\over 16}\right)^{\nu_2}\right]^m .
\label{expatalk}
\ee
The second order non linear ODE  for $u(x)$ can be recast into  a system of two integral equations and solved by successive approximations in exactly the same way as in Shimomura's works \cite{Sh}. As a result, the following theorem  holds \cite{D3}
\bth
   For any 
complex  $\nu_1$, $\nu_2$ such that $ 
\Im \nu_2\neq 0
$
there exist a  positive number $r<1$  and a  transcendent 
\be
\label{tra}
 y(x,\nu_1,\nu_2) = \wp \Bigl(\nu_1\omega_1(x)+\nu_2\omega_2(x) +v(x,\nu_1,\nu_2);~ \omega_1(x),\omega_2(x)
 \Bigr)+{1+x\over 3}.
\ee 
such that $v(x,\nu_1,\nu_2)$ is a holomorphic function  on 
 $ {\cal D}(r;\nu_1,\nu_2)$, with convergent  expansion (\ref{expatalk}). 
The coefficients $a_n$, 
$b_{nm}$ and $c_{nm}$ are certain 
rational functions of $\nu_2$. For any complex $\nu_1$ and real $\nu_2$ such that $0<\nu_2<1$, there exist a  positive number $r<1$  and a  transcendent (\ref{tra}) 
such that $v(x;\nu_1,\nu_2)$ is holomorphic with convergent  expansion (\ref{expatalk}) in 
 $ {\cal D}_0(r)$.  If $1<\nu_2<2$ the same result holds with $v=v(x;-\nu_1,2-\nu_2)$ in 
 $ {\cal D}_0(r)$. 
\label{trattoria}
\eth 
Combining the expansion (\ref{expatalk}), the expansion at $x=0$ of $\omega_1(x)$ and $\omega_2(x)$ and   the Fourier expansion of the $\wp$ function, we obtain the critical behaviours of $y(x,\nu_1,\nu_2)$ when $x\to 0$ along a regular path inside the domain ${\cal D}(r;\nu_1,\nu_2)$ or ${\cal D}_0(r)$.  
We do not rewrite the behaviours here, because they coincide with (\ref{locintro1}), (\ref{hoshi1}) and (\ref{invasi}),  with the identification 
$$
\sigma=1-\nu_2,~~~\hbox{( or } \sigma=\nu_2-1 \hbox{ when } 1<\nu_2<2 \hbox{)}.
$$
Also this result  proves - as in section \ref{shimomura} - the  convergence of the full expansion corresponding to the behaviours (\ref{locintro}), (\ref{hoshi0}) and (\ref{invos1})


\section{Isomonodromy  Deformation on the Universal Covering of a Punctured Neighborhood of a Critical Point}
 \label{isoest}

 A local analysis similar to that of section \ref{shimomura}  can be also applied  to the Schlesinger equations associated to PVI. In this way, the  result of lemma 2.4.8  of \cite{SMJ}, which yields the critical behaviours (\ref{Lambda})   and thus (\ref{locintro}) and (\ref{hoshi0}), turns out to hold  when $x\to 0$  in the domains ${\cal D}_s(r;\sigma,a)$ of section \ref{shimomura}, for   $\sigma\not\in (-\infty,0]\cup[1,+\infty)$.  The proof of this fact can be found in \cite{D4}. The critical behaviours of $A_0(x)$, $A_x(x)$ and $A_1(x)$ obtained in this way provide  (\ref{locintro1}) and (\ref{HOshi1}) by substitution into (\ref{ipsilon}),  and then (\ref{invos11}) by a transformation (\ref{sym2}).

As a consequence,  the monodromy preserving deformation approach works  for the transcendents in the Shimomura and  the Elliptic representations.    For Shimomura transcendents   the formula providing the monodromy exponent $\sigma$  at $x=0$ in (\ref{shimo}) is again 
\be
\label{rosar}
2\cos\pi \sigma=p_{0x},~~~~~\sigma\not\in (-\infty,0]\cup[1,+\infty).
\ee
In the Elliptic representation, the analogous formula is \cite{D3}
\be
\label{nudue}
 2\cos\pi\nu_2=-p_{0x},~~~~~\nu_2\not\in(-\infty,0]\cup\{1\}\cup[2,+\infty),
\ee
The constant $a$ in (\ref{shimo}) is given by a  formula which coincides with  $a_0$ in (\ref{STA1}).   In the following, we rewrite it in the simple form 
$$
a=a(\sigma).
$$

\begin{figure}
\centerline{\includegraphics[width=0.5\textwidth]{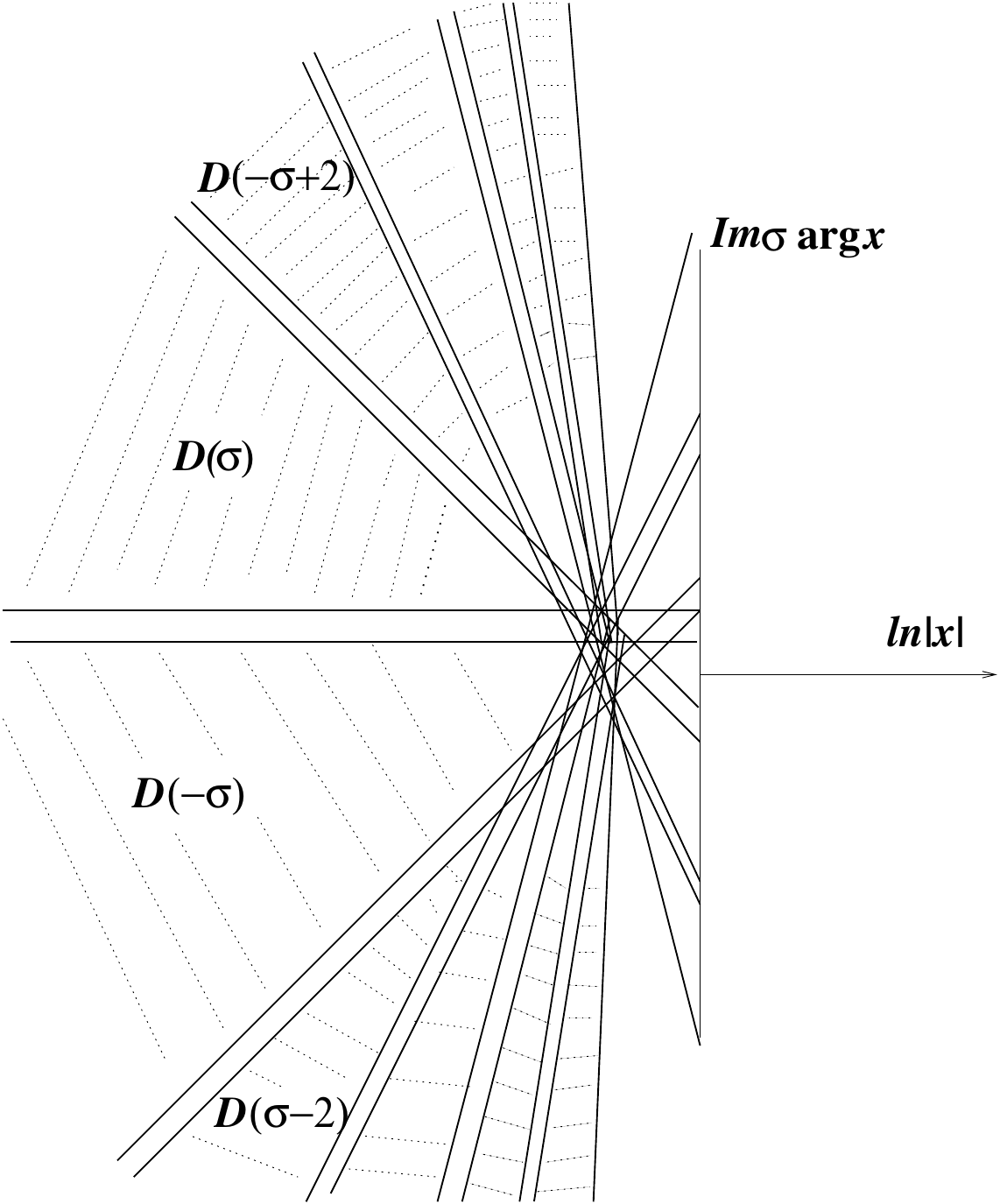}}
\caption{The union (dashed region) of some of the domains ${\cal D}_s(r_N^\pm,\sigma_N^{\pm},a(\sigma_N^{\pm}))$.}
\label{GIG5}
\end{figure}

\begin{figure}
\centerline{\includegraphics[width=0.5\textwidth]{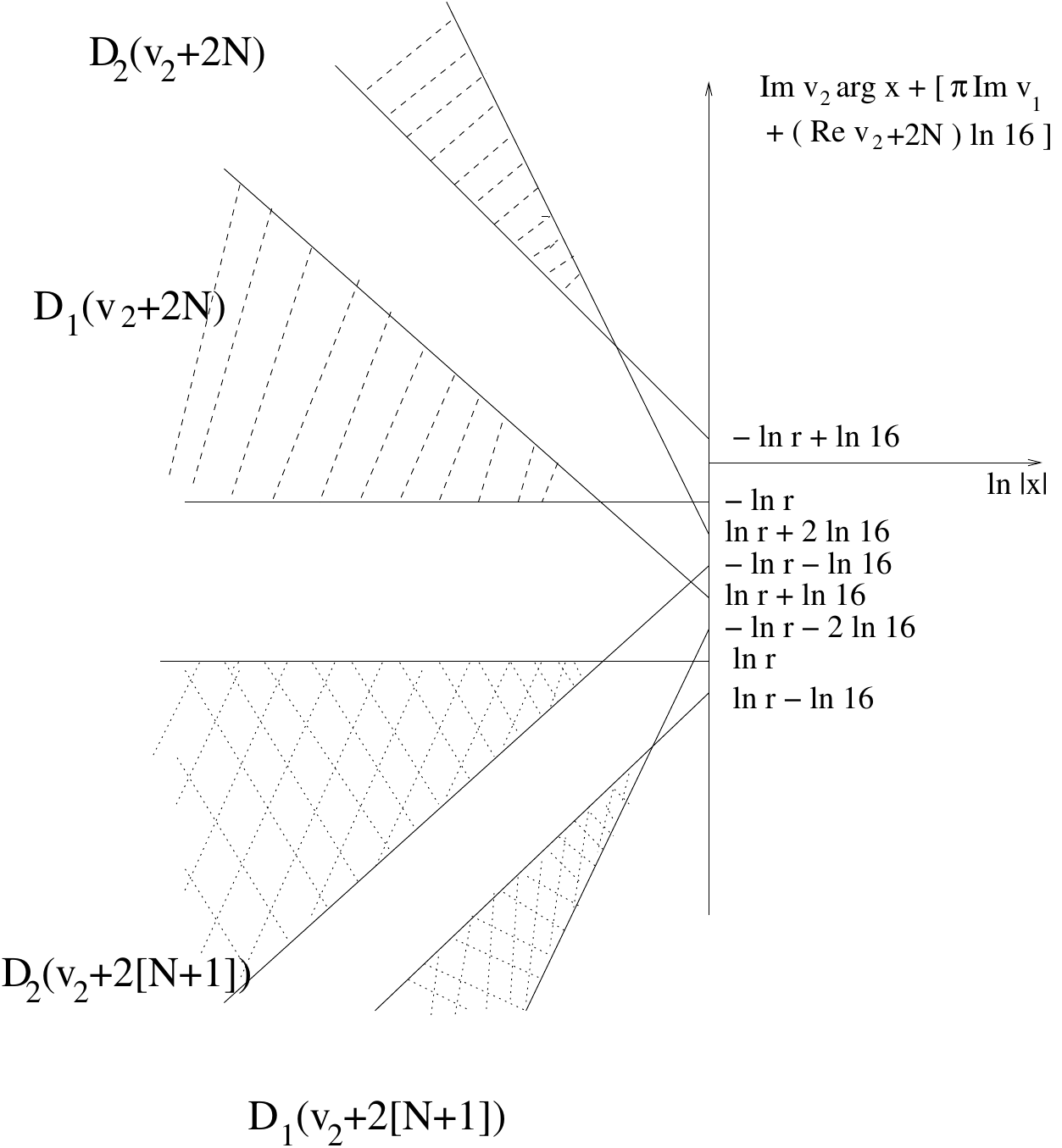}}
\caption{The domains ${\cal D}_1(r;\nu_1,\nu_2+2N):={\cal D}(r;\nu_1,\nu_2+2N)$, 
${\cal D}_2(r;\nu_1,\nu_2+2N):={\cal D}(r;-\nu_1,2-\nu_2-2N)$ 
and  ${\cal D}_1(r;\nu_1,\nu_2+2[N+1])$, 
${\cal D}_2(r;\nu_1,\nu_2+2[N+1])$ for  arbitrarily 
 fixed values of $\nu_1$, $\nu_2$, $N$. They are represented in the plane  
$(\ln|x|,~\Im \nu_2 \arg x +[\pi \Im \nu_1 +(\Re \nu_2 +2N)\ln 16])$.  A PVI transcendent can be represented  on the union  of the domains, with different $\nu_1$ and $v(x)$ on each domain.}
\label{figur1}
\end{figure}

\subsection{Extension of the Domains} 

  The isomonodromy approach on the domains ${\cal D}_s(r;\sigma,a)$  has an important consequence: for a given PVI and  given monodromy data $p_{0x},p_{x1}, p_{01}$, let us consider  (\ref{shimo}).  
If $0\leq \Re \sigma <1$, then $y(x,a,\sigma)$  coincides with a branch (\ref{locintro}) or (\ref{hoshi0}) at $x=0$ (Jimbo's result with $a=a_0$ , $\sigma=\sigma_{0x}$). If we choose another solution of (\ref{rosar}), namely 
$$ 
     \sigma_N^{\pm}:= \pm\sigma+2N, ~~~~~N\in{\bf Z}
$$
then  the branch associated to the monodromy data  $p_{0x},p_{x1}, p_{01}$  extends to the universal covering of a punctured neighborhood of $x=0$  as a PVI-transcendent with Shimomura representations 
$$
 y(x, \sigma_N^{\pm},a(\sigma_N^{\pm}))= {1 \over \cosh^2({\sigma_N^{\pm}-1\over 2}\ln {x}
    -{1\over 2}\ln a(\sigma_N^{\pm})+\ln 2 +{v({x})\over 2})},
$$
on the different domains  ${\cal D}_s(r_N;\sigma_N^{\pm},a(\sigma_N^{\pm}))$,  
with different monodromy exponents $\sigma_N^{\pm}$ and different $a=a(\sigma_N^{\pm})$. See figure \ref{GIG5}. Note that $r=r_N$.

In the Elliptic representation, the same fact is related to   the periodicity and parity of the $\wp$ function, which  allow  to change $(\nu_1,\nu_2)$ to $(\nu_1,\nu_2+2N)$ and to $(-\nu_1,2-\nu_2-2N)$, where $N\in {\bf Z}$.  
 The integration constant $\nu_2$ is expressed by (\ref{nudue}), 
which is invariant for $\nu_2\mapsto \nu_2+2N$ and  $\nu_2\to 2-\nu_2-2N$.  On the other hand   $\nu_1$ - like $a$ -   depends on the monodromy data by an explicit parametric formula (see \cite{D4} \cite{D3}) 
\be
\label{nuuno}
 \nu_1=\nu_1(\nu_2,p_{x1},p_{01}),
\ee
Now, let us consider a transcendent (\ref{tra}) of Theorem \ref{trattoria} and do the substitution $\nu_2\mapsto \nu_2+2N$ or $\nu_2\to 2-\nu_2-2N$. Provided $r=r_N$ is suitably small, the theorem still holds.  
Provided  that also  we change $\nu_1=\nu_1(\nu_2,p_{x1},p_{01})\mapsto \nu_1( \nu_2+2N,p_{x1},p_{01}) $ or $\nu_1(2-\nu_2-2N,p_{x1},p_{01}) $, we obtain different representations of the same transcendent on the  domains
  ${\cal D}(r_N;\nu_1,\nu_2+2N)$ 
 and  the domains ${\cal D}(r_N,-\nu_1,2-\nu_2-2N)$. 
%
%
We  draw the domains in the $(\ln|x|,\Im\nu_2 \arg x)$-plane, in  figure \ref{figur1}.

\subsection{The Final Picture on the Universal Covering of a Punctured Neighborhood of a Critical Point}

\begin{figure}
\centerline{\includegraphics[width=0.5\textwidth]{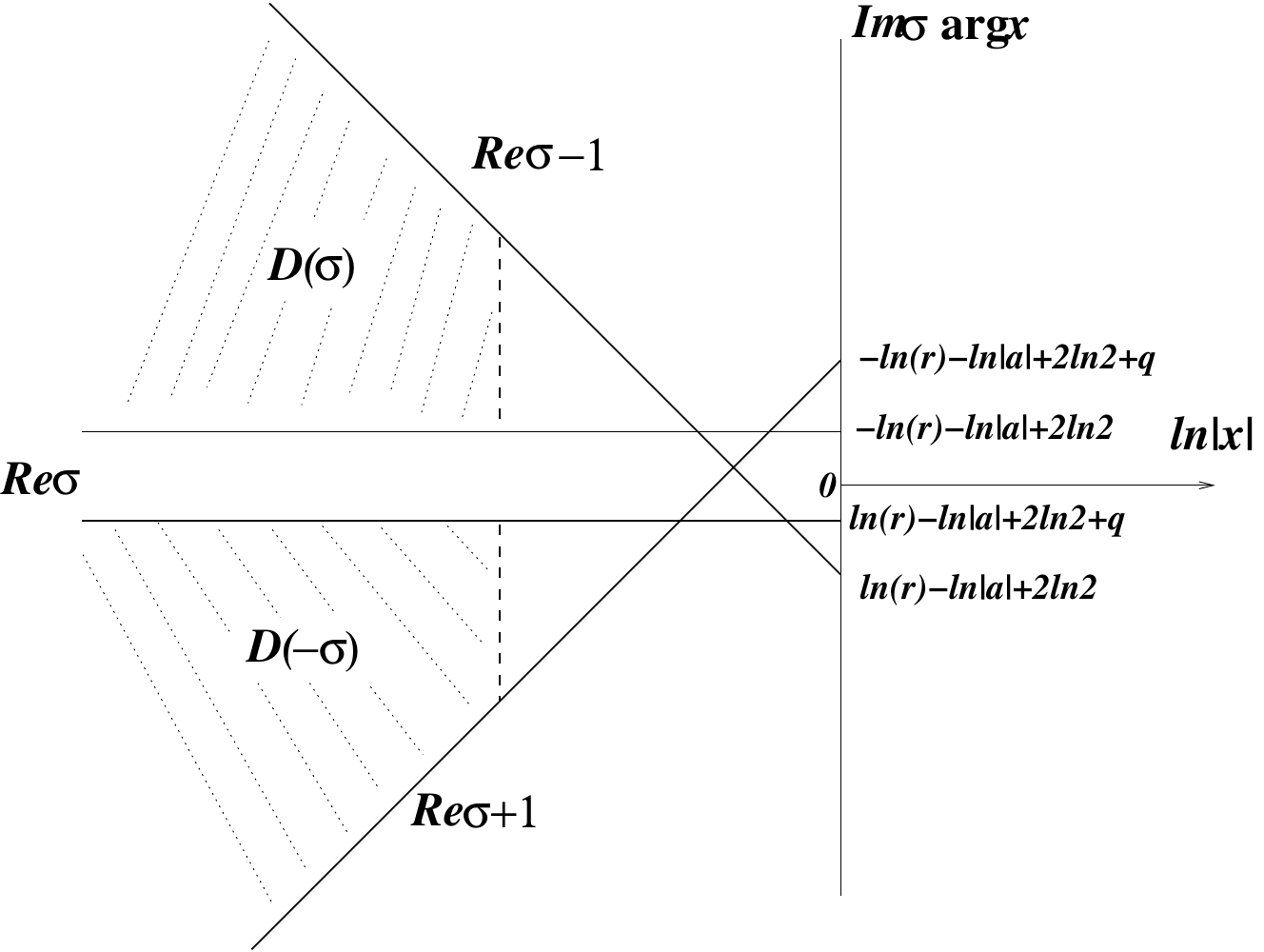}}
\caption{The domains ${\cal D}_s(r,\sigma,a(\sigma))$ and ${\cal D}_s(r;-\sigma,a(-\sigma))$, in the $(\ln |x|, \Im \sigma \arg x)$ plane. The slopes of the boundaries are indicated. In the particular case in figure, $\Re \sigma=0$. $q=\ln 
{|\sigma^2| \over 4|A|^2}$.}
\label{GIG1}
\end{figure}

\begin{figure}
\centerline{\includegraphics[width=0.5\textwidth]{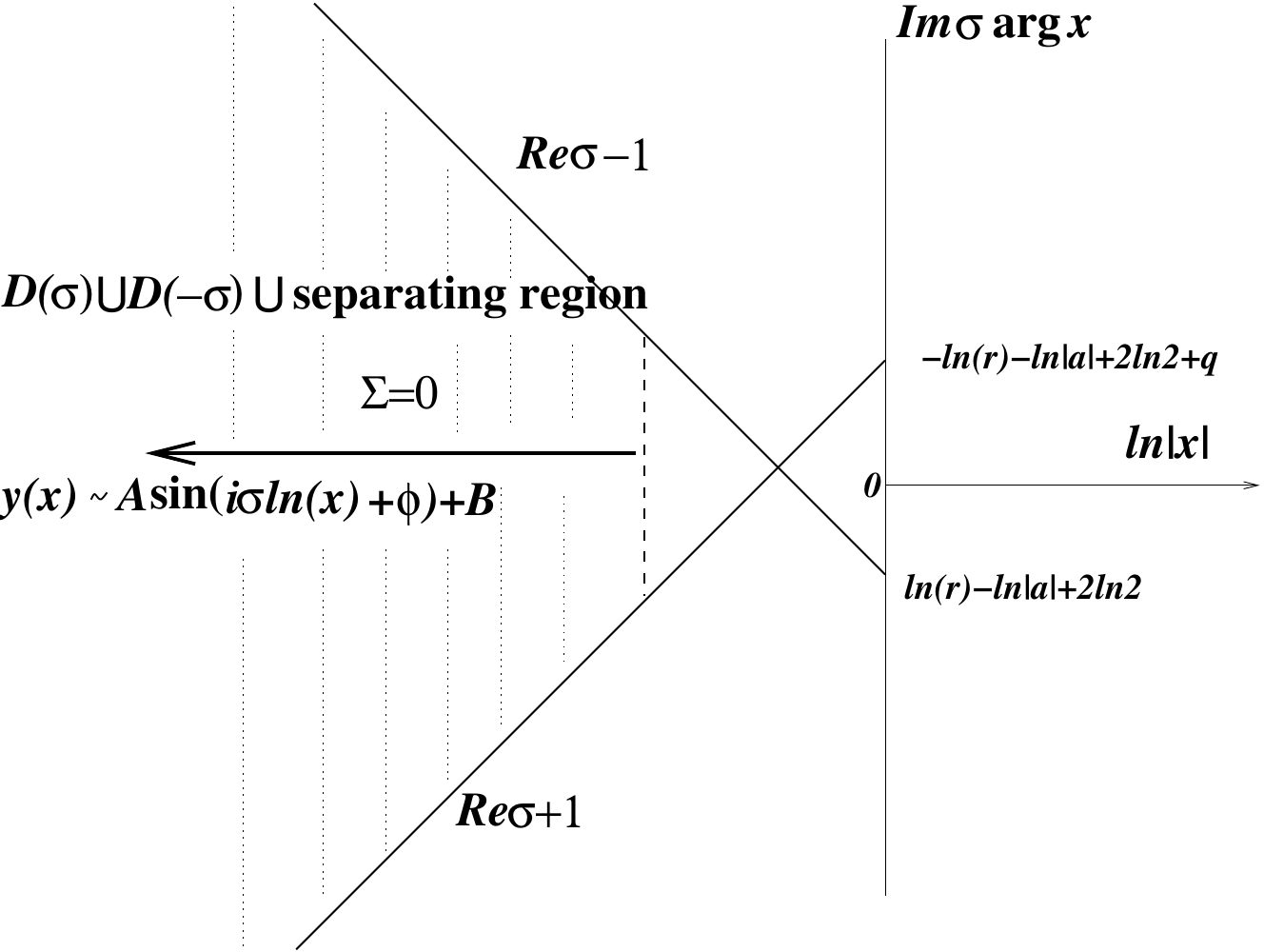}}
\caption{The union ${\cal D}_s(r;\sigma,a(\sigma))\cup\hbox{\{separating region\}}\cup {\cal D}_s(r;-\sigma,a(-\sigma))$. Along the  paths with $\Sigma=0$  the behaviour of $y(x)$ is (\ref{be1})=(\ref{be2})=(\ref{HOshi1}). This figure corresponds to the case of  figure \ref{GIG1}.}
\label{GIG2}
\end{figure}

\begin{figure}
\centerline{\includegraphics[width=0.5\textwidth]{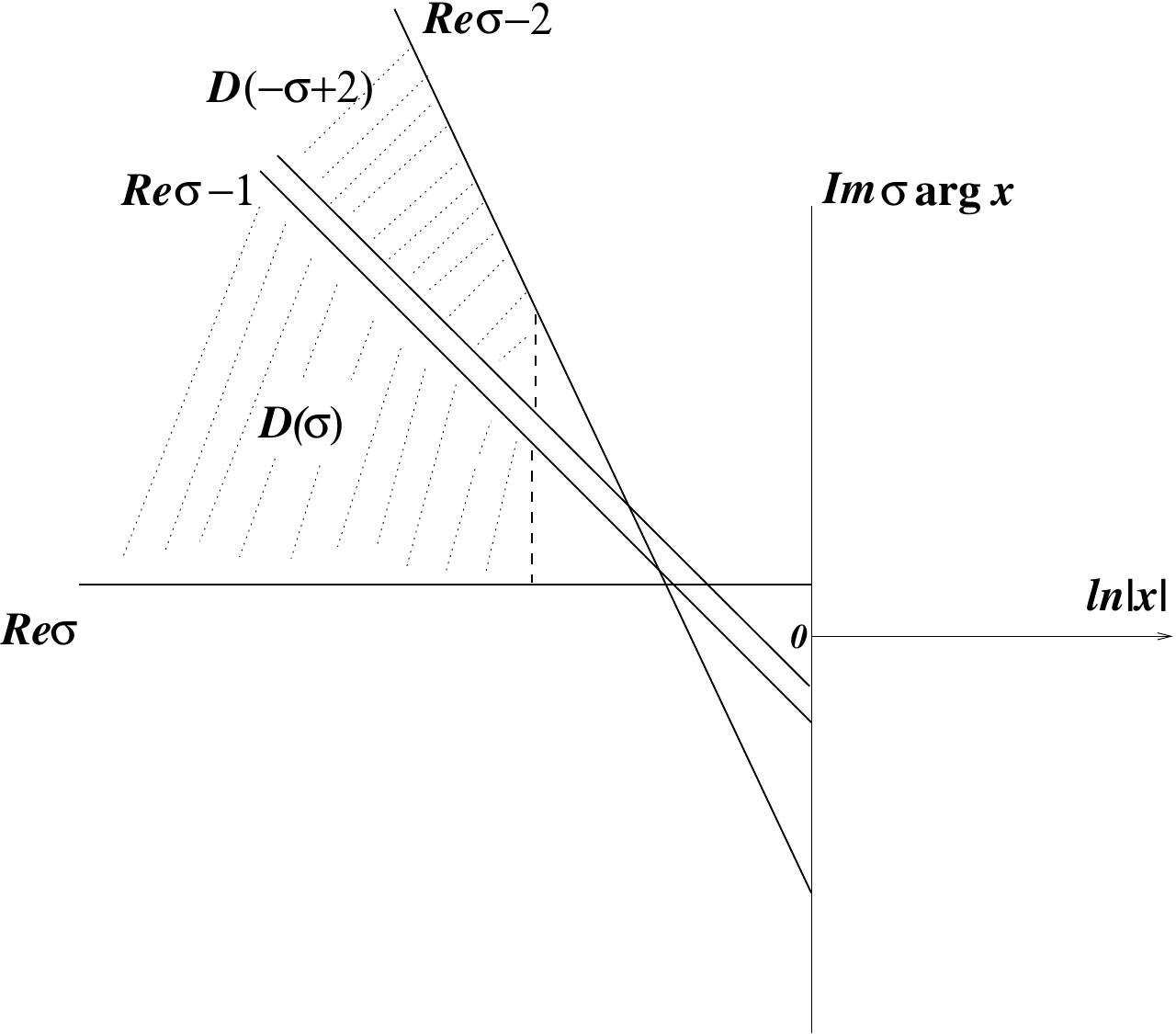}}
\caption{The domains ${\cal D}_s(r,\sigma,a(\sigma))$ and ${\cal D}_s(r,-\sigma+2,a(-\sigma+2))$. The slopes of the boundaries are indicated. This figure corresponds to the case of  figure \ref{GIG1}.}
\label{GIG3}
\end{figure}

\begin{figure}
\centerline{\includegraphics[width=0.5\textwidth]{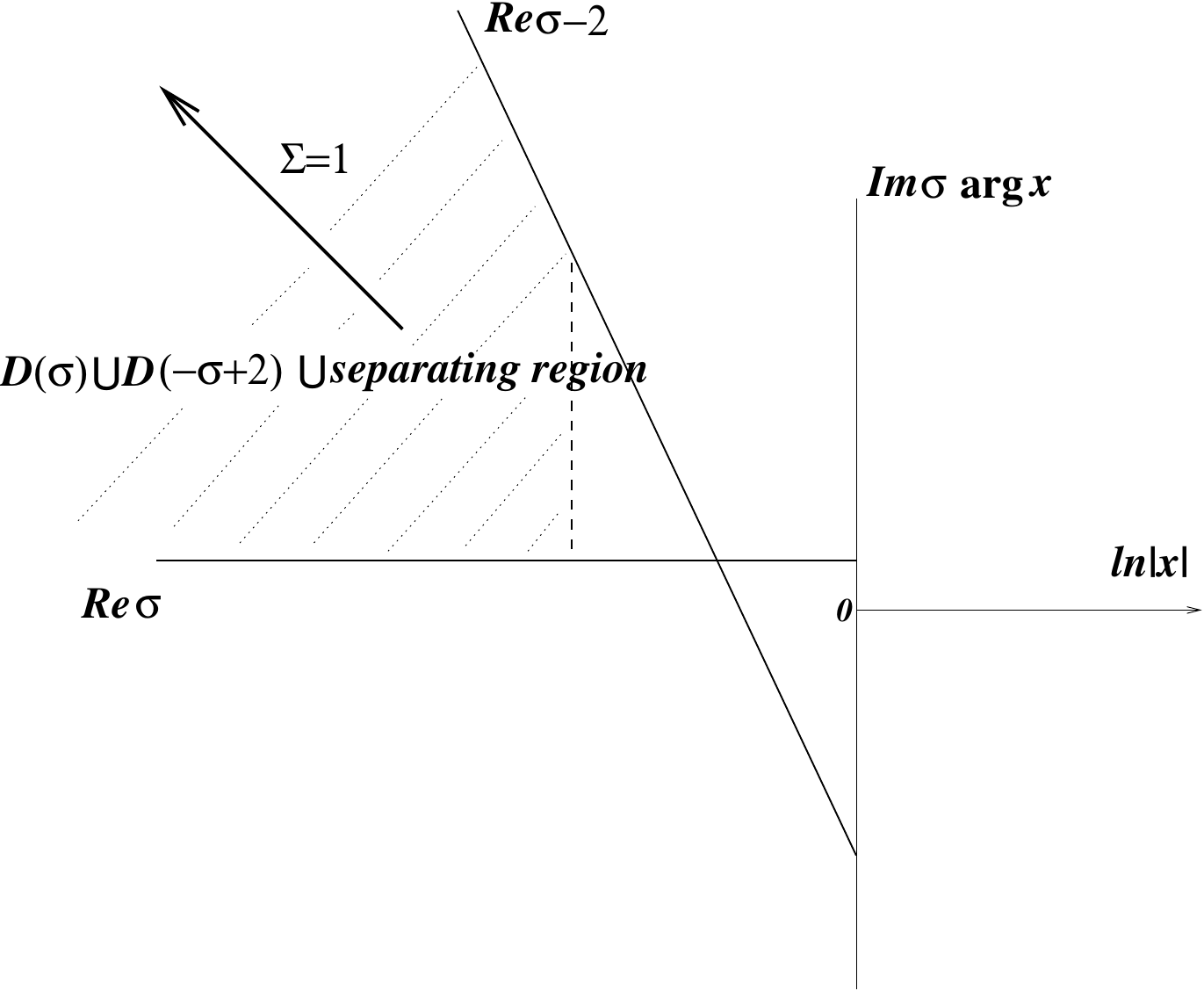}}
\caption{The union ${\cal D}_s(r,\sigma,a(\sigma))\cup\hbox{\{separating region\}}\cup{\cal D}_s(r,-\sigma+2,a(-\sigma+2))$, and the paths for  $\Sigma=1$, along which the behaviour of $y(x)$ is (\ref{invos11}). This figure corresponds to the case of  figure \ref{GIG3}.}
\label{GIG4}
\end{figure}

\begin{figure}
\centerline{\includegraphics[width=0.5\textwidth]{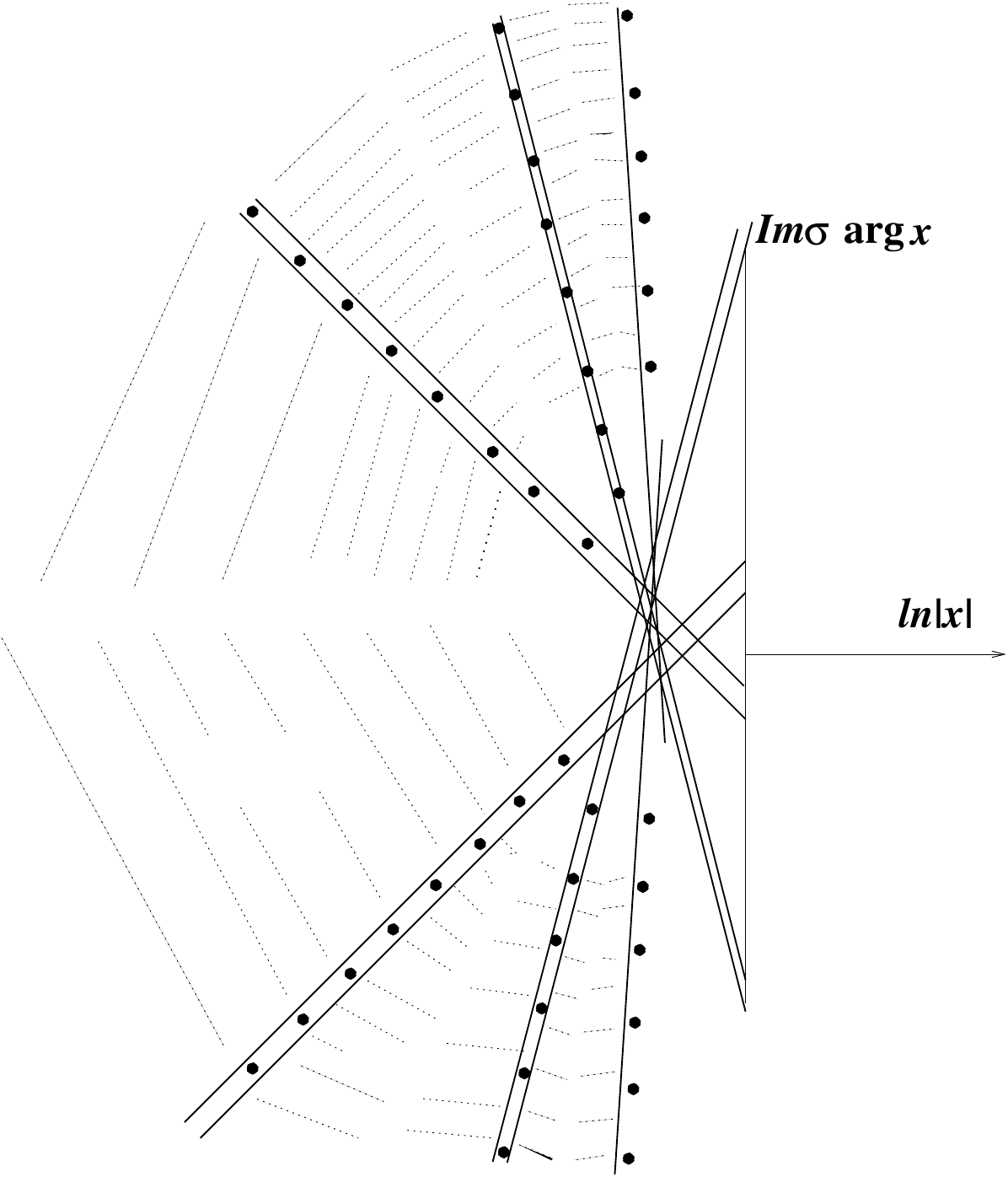}}
\caption{The possible location of the poles (indicated by dots). This figure refers to the case of figure \ref{GIG5}. Note that in this particular figure $\Re \sigma=0$, therefore there are definitively no poles if $x\to 0$ with $\arg(x)$ bounded, which instead exist when $\Re\sigma=1$.}
\label{GIG6}
\end{figure}

We  can further improve   the description of the critical behaviours of a transcendent on the universal covering of the punctured neighbourhood of $x=0$.  Let us study what happens when we pass from ${\cal D}_s(r;\sigma,a(\sigma))$ to its neighbouring domains ${\cal D}_s(r;-\sigma,a(-\sigma))$ and ${\cal D}_s(r,-\sigma+2,a(-\sigma+2))$. 

To start,  consider the two Shimomura representations   $y(x,\sigma,a(\sigma))$ and $y(x,-\sigma,a(-\sigma))$ of a given transcendent, defined on neighbouring domains ${\cal D}_s(r;\sigma,a(\sigma))$  and  ${\cal D}_s(r;-\sigma,a(-\sigma))$ respectively. From the parametric connection formulae for $a=a(\sigma)$  one can prove that   
\be
\label{Aya}
a(\sigma) ~a(-\sigma)= {4\sigma^2\over A^2(\sigma^2)}
\ee
where $A^2(\sigma^2)$ is (\ref{AB}).  The above implies that 
 ${\cal D}_s(r;-\sigma,a(-\sigma))$ is 
$$
(\Re \sigma +1) \ln |x|+\{2\ln2-\ln|a|+\ln{|\sigma^2|\over 4|A|^2}\}-\ln r<\Im \sigma \arg x < ~~~~~~~~~~
$$
$$~~~~~~~~~~~~~<\Re \sigma\ln |x|+\{2\ln 2-\ln|a|+\ln {|\sigma^2|\over 4|A|^2}\}+\ln r,
$$
while, by definition, 
${\cal D}_s(r;\sigma,a(\sigma))$ is 
$$
\Re \sigma \ln |x|+\{2\ln2-\ln|a|\}-\ln r<\Im \sigma \arg x < (\Re \sigma-1)\ln |x|+\{2\ln 2-\ln|a|\}+\ln r,
$$
The two domains  do not in general intersect, as in figure \ref{GIG1}.   The critical behaviours inside the domains, along paths with $\Sigma=0$,  are of type (\ref{HOshi1}): 
\be
\label{be1}
  y(x,\sigma,a(\sigma))=x\{ A^+(\sigma^2)\sin(i\sigma\ln x +\phi(\sigma))+B(\sigma^2)\}+O(x^2)
\ee 
when $x_0\in {\cal D}_s(r,\sigma,a(\sigma))$, and 
\be
\label{be2}
y(x,-\sigma,a(-\sigma))=x\{ A^-(\sigma^2)\sin(-i\sigma\ln x +\phi(-\sigma))+B(\sigma^2)\}+O(x^2)
\ee
when $x_0\in{\cal D}_s(r;-\sigma,a(-\sigma))$. Here 
 $A^\pm(\sigma^2)$ are the roots of $A^2(\sigma^2)$ in (\ref{AB}), and  $B(\sigma^2)$ too is given in (\ref{AB}). Also,  $a(\sigma) =Ae^{i\phi(\sigma)}/2i$. Due to the arbitrary sign we can take $A^-(\sigma^2)=-A^+(\sigma^2)$, so that formula (\ref{Aya}) implies $\phi(-\sigma)=-\phi(\sigma)$.  It follows that    (\ref{be1}) and  (\ref{be2}) are formally the same when $\Sigma=0$, though on different domains. 
On the other hand, if  $x_0$ belongs to the region  separating the two domains, the  only path  (\ref{PATH}) which does not  enter into one of the domains when $x\to 0$,  is precisely a path with  $\Sigma=0$.
In \cite{D4} we proved  that  {\it (\ref{be1}) and  (\ref{be2})  hold also in the separating region}.   Thus, 
{\it  the behaviour (\ref{HOshi1}) extends to  ${\cal D}_s(r;\sigma,a(\sigma))\cup\hbox{\{separating region\}}\cup {\cal D}_s(r;-\sigma,a(-\sigma))$}. 
 This is represented in figure \ref{GIG2}. 

\vskip 0.2 cm 

\noindent
{\bf Remark:} 
The fact that  domains ${\cal D}_s(r;\sigma,a(\sigma))$ and ${\cal D}_s(r;-\sigma,a(-\sigma))$  do not in general intersect means that the series  (\ref{typeseries})   of $v(s,\sigma,a(\sigma))$ and $v(x,-\sigma,a(-\sigma))$ do not converge on the region that separates ${\cal D}_s(r;\sigma,a(\sigma))$ and ${\cal D}_s(r;-\sigma,a(-\sigma))$.

\vskip 0.2 cm 
The above construction can be repeated for the neighbouring domains ${\cal D}_s(r,\sigma,a(\sigma))$ and ${\cal D}_s(r,-\sigma+2,a(-\sigma+2))$ of figure \ref{GIG3}, making use of the symmetry (\ref{sym2}) applied to (\ref{HOshi1}). It   transforms  $p_{0x}\mapsto -p_{0x}$, therefore  $\sigma\mapsto 1-\sigma$ and (\ref{HOshi1}) into (\ref{invos11}), to be understood as defined 
on the union    ${\cal D}_s(r,\sigma,a(\sigma))\cup\hbox{\{separating region\}}\cup{\cal D}_s(r,-\sigma+2,a(-\sigma+2))$, as in figure \ref{GIG4}. 
The denominator does not vanish in ${\cal D}_s(r,\sigma,a(\sigma))$ and ${\cal D}_s(r,-\sigma+2,a(-\sigma+2))$, where it has convergent expansion.  It may vanish in the separating region. Thus, the poles of (\ref{invos11}) possibly lie  in the separating region. The zeros of the leading term ${\cal A}\sin(i(1-\sigma)\ln x+\phi)+{\cal B}$ are computed in a form similar to  (\ref{seRO}).   If we write the {\it formal} expansion of the denominator  we formally  compute the equivalent of  (\ref{POLESintro}).

 \vskip 0.2 cm 
If we consider the union  in figure \ref{GIG5} of the domains ${\cal D}_s(r_N^\pm,\sigma_N^{\pm},a(\sigma_N^{\pm}))$ (up to a fixed ``big'' $N$) and apply the above considerations to neighbouring domains, we obtain the final picture:  a  transcendent  defined on the universal covering of a punctured neighbourhood of $x=0$ behaves as follows:

When $x\to 0$ along lines of slope $\Re\sigma_N^+$ in the $(\ln |x|,\Im\sigma\arg x)$ plane, the behaviour is of type (\ref{HOshi1}) , with $\sigma\mapsto \sigma_{N}^{+}$.

 When $x\to 0$ along lies with slope $\Re\sigma_N^+ -1$ the behaviour is of type (\ref{invos11}), with $\sigma\mapsto \sigma_N^+$. 

When $x\to 0$ with other slopes, the behaviour is of type (\ref{locintro1}), with $\sigma\mapsto \sigma_N^\pm$. 

\noindent
Asymptotically, the  poles  possibly lie along lines with slopes $\Re\sigma+2l-1$, $l\in {\bf Z}$,  in the regions separating the domains, as in figure \ref{GIG6}. In the particular case  $\Re\sigma=1$ and  $l=0$, the poles are described in   section \ref{Poles}.

The results at $x=0$ are transferred to $x=1$ and $x=\infty$ by the symmetries (\ref{onara}) and (\ref{onara1}).


\end{document}